\documentclass{amsart}
\usepackage[dvips]{graphicx}
\usepackage{amscd}
\usepackage{amsmath}
\usepackage{amsxtra}
\usepackage{amsfonts}
\usepackage{amssymb}

\newtheorem{theorem}{Theorem}[section]
\newtheorem{corollary}[theorem]{Corollary}
\newtheorem{lemma}[theorem]{Lemma}
\newtheorem{proposition}[theorem]{Proposition}

\theoremstyle{definition}
\newtheorem{definition}[theorem]{Definition}
\newtheorem{remark}[theorem]{Remark}

\newtheorem{example}[theorem]{Example}
\theoremstyle{remark}

\renewcommand{\theclaim}{\textup{\theclaim}}

\newtheorem*{acknowledgements}{Acknowledgements}

\numberwithin{equation}{section}

\def\openone

{\mathchoice

{\hbox{\upshape \small1\kern-3.3pt\normalsize1}}

{\hbox{\upshape \small1\kern-3.3pt\normalsize1}}

{\hbox{\upshape \tiny1\kern-2.3pt\SMALL1}}

{\hbox{\upshape \Tiny1\kern-2pt\tiny1}}}

\makeatletter

\newbox\ipbox

\newcommand{\ip}[2]{\left\langle #1\, , \,#2\right\rangle}
\newcommand{\diracb}[1]{\left\langle #1\mathrel{\mathchoice

{\setbox\ipbox=\hbox{$\displaystyle \left\langle\mathstrut
#1\right.$}

\vrule height\ht\ipbox width0.25pt depth\dp\ipbox}

{\setbox\ipbox=\hbox{$\textstyle \left\langle\mathstrut
#1\right.$}

\vrule height\ht\ipbox width0.25pt depth\dp\ipbox}

{\setbox\ipbox=\hbox{$\scriptstyle \left\langle\mathstrut
#1\right.$}

\vrule height\ht\ipbox width0.25pt depth\dp\ipbox}

{\setbox\ipbox=\hbox{$\scriptscriptstyle \left\langle\mathstrut
#1\right.$}

\vrule height\ht\ipbox width0.25pt depth\dp\ipbox}

}\right. }

\newcommand{\dirack}[1]{\left. \mathrel{\mathchoice

{\setbox\ipbox=\hbox{$\displaystyle \left.\mathstrut
#1\right\rangle$}

\vrule height\ht\ipbox width0.25pt depth\dp\ipbox}

{\setbox\ipbox=\hbox{$\textstyle \left.\mathstrut
#1\right\rangle$}

\vrule height\ht\ipbox width0.25pt depth\dp\ipbox}

{\setbox\ipbox=\hbox{$\scriptstyle \left.\mathstrut
#1\right\rangle$}

\vrule height\ht\ipbox width0.25pt depth\dp\ipbox}

{\setbox\ipbox=\hbox{$\scriptscriptstyle \left.\mathstrut
#1\right\rangle$}

\vrule height\ht\ipbox width0.25pt depth\dp\ipbox}

} #1\right\rangle}

\newcommand{\bz}{\mathbb{Z}}

\newcommand{\br}{\mathbb{R}}

\newcommand{\bt}{\mathbb{T}}
\newcommand{\bn}{\mathbb{N}}

\newcommand{\vectr}[2]{\left[\begin{array}{r}#1\\#2\end{array}\right]}

\def\blfootnote{\xdef\@thefnmark{}\@footnotetext}

\renewcommand{\mod}{\operatorname{mod}}

\input xy
\xyoption{all}
\usepackage{amssymb}

\def\-{^{-1}}

\def\sa{\mathcal S_A}

\def\rc{\mathcal R_C}

\def\at{(A^T)}
\def\bza{\bz^d[A^{-1}]}
\newcommand{\dist}{\operatorname*{dist}}

\begin{document}
\title[Wavelet groups]{UNITARY REPRESENTATIONS OF WAVELET GROUPS AND ENCODING OF ITERATED FUNCTION SYSTEMS IN SOLENOIDS.}
\author{Dorin Ervin Dutkay}
\blfootnote{Research supported in part by a grant from the National Science Foundation DMS-0704191}
\address{[Dorin Ervin Dutkay]University of Central Florida\\
	Department of Mathematics\\
	4000 Central Florida Blvd.\\
	P.O. Box 161364\\
	Orlando, FL 32816-1364\\
U.S.A.\\} \email{ddutkay@mail.ucf.edu}
\author{Palle E.T. Jorgensen}
\address{[Palle E.T. Jorgensen]University of Iowa\\
        Department of Mathematics\\
        14 MacLean Hall\\
        Iowa City, IA, 52242-1419\\
U.S.A.\\}
\email{jorgen@math.uiowa.edu}
\author{Gabriel Picioroaga}
\address{[Gabriel Picioroaga]Binghamton University\\
         Department of Mathematical Sciences\\
         Binghamton, NY, 13902-6000\\
U.S.A.\\} 
\email{gabriel@math.binghamton.edu}
\thanks{} 
\begin{abstract}
For points in $d$ real dimensions, we introduce a geometry for general digit sets. We introduce a positional number system where the basis for our representation is a fixed $d$ by $d$ matrix over $\bz$. Our starting point is a given pair $(A, \mathcal D)$ with the matrix $A$ assumed expansive, and $\mathcal D$ a chosen complete digit set, i.e., in bijective correspondence with the points in $\bz^d/A^T\bz^d$. We give an explicit geometric representation and encoding with infinite words in letters from $\mathcal D$. We show that the attractor $X(A^T,\mathcal D)$ for an affine Iterated Function System (IFS) based on $(A,\mathcal D)$ is a set of fractions for our digital representation of points in $\br^d$. Moreover our positional ``number representation'' is spelled out in the form of an explicit IFS-encoding of a compact solenoid $\sa$ associated with the pair $(A,\mathcal D)$.  The intricate part (Theorem \ref{thenccycl}) is played by the cycles in $\bz^d$ for the initial $(A,\mathcal D)$-IFS. Using these cycles we are able to write down formulas for the two maps which do the encoding as well as the decoding in our positional $\mathcal D$-representation.

We show how some wavelet representations can be realized on the solenoid, and on symbolic spaces. 
\end{abstract}
\maketitle \tableofcontents

\section{Introduction}

Let $A$ be a $d\times d$ matrix over $\bz$ (the integers), and assume that its eigenvalues $\lambda$ satisfy $|\lambda|>1$. In particular $A$ is assumed invertible. Let $\bz^d[A^{-1}]$
be the associated discrete group obtained as an inductive limit
 $$\bz^d\stackrel{A}{\longrightarrow}\bz^d\stackrel{A}{\longrightarrow}\bz^d\stackrel{A}{\longrightarrow}\dots$$
 i.e., $\bz^d[A^{-1}]=\cup_{k=0}^\infty A^{-k}\bz^d$, 
  and let $\mathcal S_A$ denote the corresponding dual compact abelian group, where duality is in the sense of Pontryagin; $\mathcal S_A$ is a {\it solenoid}.

  Let $A^T$ be the transposed matrix, and let $i:\bz^d[A^{-1}]\rightarrow\br^d$ be the natural embedding (of groups), so $i$ is a homomorphism of the respective additive groups. Now let $\hat i:\br^d\rightarrow\mathcal S_A$ be the dual homomorphism. Hence $\hat i$  embeds $\br^d$ as a subgroup of the compact solenoid $\mathcal S_A$. We will need this generality in an application to the analysis of wavelet multiresolutions. Here we use $\widehat{\br^d}=\br^d$, i.e., $\br^d$ is its own Pontryagin dual. In a special case of this construction, for $d = 1$, our embedding $\hat i$ corresponds to ergodic theoretic flows on compact spaces studied in a variety of contexts in dynamics.
  
  Note $\hat i(A^Tx)=\sigma_A(\hat i(x))$, $x\in\br^d$, where $\sigma_A:\mathcal S_A\rightarrow\mathcal S_A$ is the endomorphism induced by $A$. 
  
    Motivated by classical number theoretic problems in digital representations of fractions, we explore here symbolic representations of points in $\mathcal S_A$. Our analysis uses an extension of George Mackey's semidirect product construction from representation theory (section 3) combined with a study of a family of combinatorial cycles (section 4.1). In section 6 we show that our encoding theorems (sections 4 and 5) apply to the construction of wavelet multiresolutions, i.e., for generalized wavelet constructions where scaling in $\br^d$ corresponds to the matrix multiplication $x\mapsto Ax$, $x\in\br^d$; and the corresponding $\bz$-action $\bz\times\br^d\ni(k,x)\mapsto A^kx\in\br^d$.

{\bf Motivation from physics: The renormalization group.} The idea of scale invariance is old. Its best know modern formulation in mathematics takes the form of iterated function systems (IFS), e.g., \cite{Hut81}; and in physics it takes the form of a Renormalization group (RG). But scaling arguments are commonplace in pure (e.g., wavelets) and applied mathematics; for example in attempts at explaining turbulence.

The renormalization group makes its appearance in physics in different guises, often as a mathematical trick to get rid of the infinities, for example in quantum field theory, see e.g., \cite{Fe87, Bal91}. As a pure technique, it obtained maturity with for example the papers \cite{Fe87, Fe94, KW88, KW91} among others. The technique was developed also in quantum electrodynamics by R. Feynman, K. Wilson and others. The physicists devised theories of mass and charge renormalization.

Old-style renormalization group (RG) techniques in physics have run into difficulties with non-renormalizability of gravity. Still they are used in various guises as tools in solid state physics, as they often get around divergence difficulties with the use of perturbation theory.

As with iterated function systems (IFS), renormalization groups in physics attempt to describe infinite systems in terms of block variables, i.e.: some magnitudes which describe the average behavior of the various constituent blocks, often approximately true in practice, and good enough, to a first approximation. This is more or less amounts to finding long term behavior of a suitable RG transformation. When iterated many times, this RG transformation leads to a certain number of fixed points analogous to those seen in (non-contractive) IFSs.

  It was suggested in the 1970s by Don Knuth and others that there is an intriguing geometry behind computations in a positional number system. As is known since Euclid, representation of numbers in a fixed basis entails expansions in powers of the chosen base, say $b$. There is then a subset $\mathcal D$ of the integers $\bz$ of cardinality $|b|$ such that the corresponding ``digital'' expansion of real numbers is encoded by finite or infinite words in the ``alphabet'' $\mathcal D$. In fact Knuth \cite{Knu76} stresses that for a fixed $b$, there are many choices of digit sets $\mathcal D$ which yield a positional number system in a sense which is made precise. Similarly Knuth suggested the use of a matrix in place of $b$. In sections 5 and 6 below, we outline both the initial suggestion, the relevant literature, starting with \cite{Odl78}; and we present our main results. The first sections of our paper address tools from representation theory central to our problem. 
  
  Both the choice of the base $b$, and the set of digits allows for a great deal of freedom, even if we restrict numbers on the real line.
           Don Knuth suggested that this idea works in higher dimensions, i.e., in encoding points in $\br^d$ this way. The case $d = 2$ of course includes positional representations of complex numbers, and associated computer generated images in the plane. But by increasing the dimension further, this suggests using, for scaling basis, instead a fixed $d$ by $d$ expansive matrix over $\bz$, say $A$, in place of the base number $b$, and a subset $\mathcal D$ of $\bz^d$ of cardinality $| \det A|$  for digits.
        This leads to the puzzling question: ``In this geometric formulation, what are then the fractions as subsets of $\br^d$ ?'' Early computer calculations by Knuth in 2D suggested that ``fractions'' take the form of ``dragon-like'' compact sets, ``Twin Dragons'' etc. Knuth's idea was taken up in later papers (by other authors, e.g., Lagarias and Wang)  under the name affine Iterated Function Systems (IFSs), and the set of fractions associated to a fixed pair $(A, \mathcal D)$ in $d$ real dimensions have been made precise in the form of attractors for the IFS defined from $(A, \mathcal D)$, see e.g., \cite{Hut81} and \cite{BrJo99}.

        For points in Euclidian space, we introduce matrix scaling and digit
sets. Our aim is to study the interplay between associated spectra and
geometry. For a fixed matrix scaling and digit set we introduce a positional
number system where the basis for our representation is a fixed $d$ by $d$
matrix $A$ over $\bz$. Specifically, a pair $(A,\mathcal D)$ is given, the matrix $A$ assumed
expansive, and a finite set $\mathcal D$ a chosen as a complete digit set, i.e., the
points in $\mathcal D$ are in bijective correspondence with the finite group $\bz^d/A\bz^d$.

    Such higher dimensional ``number systems'' allow more flexibility than the
classical one, introducing a computational device for the study of, for
example, exotic tilings, wavelet sets and fractals. These are geometric
structures in $\br^d$, studied recently in \cite{BMM99, Cho07, JKS07, Shu03, Hut81,
BHS05, FMM06, Rud89}. We take advantage of a natural embedding of $\br^d$ in an
associated solenoid, and we obtain an explicit solenoid-encoding of
geometries in $\br^d$, giving insight into notions of redundancy, and offering a
computational tool.

   This expanded view also introduces novelties such as
non-commutativity into encoding. We give an explicit geometric
representation and encoding for pairs pair $(A,\mathcal D)$ (Theorem \ref{prop3_9}), i.e., an
encoding with specific infinite words in letters from $\mathcal D$.
    
    Our positional ``number representation'' takes the form of an explicit
IFS-encoding of points in a compact solenoid SA associated with the pair
$(A,\mathcal D)$. A crucial part (Theorem \ref{thenccycl}) is played by certain extreme cycles in
the integer lattice $\bz^d$ for the initial $(A,\mathcal D)$-IFS. Using these cycles we
write down a formula for the two maps which do the encoding as well as the
decoding in our positional $\mathcal D$-representation.

        We will need basic tools from spectral theory (e.g., \cite{Arv02}), but our aim
is computational, still using operator algebraic tools, analogous to those
used in the analysis of graphs and generalized multiresolutions, see e.g.,
\cite{BMM99, Cho07, JKS07, Shu03}. Our use of iterated function systems follows
conventions from \cite{Hut81, BHS05, FMM06}.

       When an invertible matrix $A$ and a finite subset $\mathcal D$ in $\br^d$ are given, then we consider an associated finite set of affine mappings $\tau_d$, indexed by points $d$ in $\mathcal D$, as follows, $$\tau_d(x):= A^{-1}(x + d).$$ Under suitable conditions on the pair $(A,\mathcal D)$, see Definition \ref{ifs}, repeated iterations of the combined system $(\tau_d)_{d \in\mathcal D}$, then yields certain limit concepts. They take the precise form of either certain compact subsets in $\br^d$ (attractors) or the form of limiting measures, so called equilibrium measures. In this formulation, the theory was made precise by J. Hutchinson in the paper \cite{Hut81}, and this gave rise to what is now known as affine iterated function systems (IFSs). As is known, e.g.,\cite{Hut81}, to make the limit notions precise, one introduces metrics on families of compact subsets in $\br^d$, or on families of probability measures. (In their primitive form these metrics generalize the known Hausdorff distance.)
       
       In harmonic analysis a subclass of the IFSs have been studied extensively by R. Strichartz and his co-authors, see \cite{Str05, Str06}. The emphasis there is on discrete potential theory, while our present focus is on tiling and coding questions. However much of our motivation derives from harmonic analysis.

         It was also realized that wavelet algorithms can be put into this framework, and for fixed $(A, \mathcal D)$, one is led to ask for Haar wavelets, and to wavelet sets.
         In this paper we show that there is a representation theoretic framework for these constructions involving a certain discrete group that is studied in algebra under the name of the Baumslag-Solitar group.

        In addition to our identifying wavelet sets in a solenoid
encoding (section 4), in section 6 we further show that the positional
number representation for $\br^d$ associated to a given pair $(A,\mathcal D)$ takes
an algorithmic form involving both ``fractions'' and ``integer points''.
Here the fractions are represented by a Hutchinson attractor $X(A^T,\mathcal D)$
and the ``$(A,\mathcal D)$-integers'' by a certain lattice $\Gamma$ which makes
$X(A^T,\mathcal D)$ tile $\br^d$ by $\Gamma$ translations. To compute this lattice $\Gamma$
which makes an $X(A^T,\mathcal D)$ tiling, we use a certain spectral duality (Lemma \ref{lemfuglede}).
    Finally (section 6.2) the last mentioned duality is illustrated
with specific planar examples, where both cycles and lattices are
worked out.

         \section{Iterated Function Systems}
 The setting of this paper is a fixed $d$ by $d$ matrix $A$ over the integers $\bz$, satisfying a certain spectral condition, and its relation to the rank-$d$ lattice $\bz^d$ in $\br^d$.

Traditional wavelet bases in $L^2(\br^d)$ are generated by a distinguished finite family of functions in $L^2(\br^d)$ and the operations translation with the rank-$d$ lattice $\bz^d$, and scaling with powers of $A^j$, $j \in \bz$. But there are similar wavelet constructions, super wavelets, in other Hilbert spaces which we explore here. See also \cite{BJMP05,DuJo06,DuJo07}

From the initial matrix $A$ we form the discrete group $\bza$ generated by powers of $A^j$, $j\in\bz$, applied to $\bz^d$; and the compact dual solenoid group $\sa$. This solenoid is related to $[0,1)^d \times \Omega$  where $\Omega$ is a compact infinite product of a fixed finite alphabet. But the two $\sa$ and $[0,1)^d \times \Omega$  are different, and their relationships are explored below.

        First recall that matrix multiplication by $A$ induces an automorphism $\sigma_A$ in $\sa$. We are interested in periodic points for this action, and in a certain family of extreme orbits called cycles. As is well known, $\br^d$ is naturally embedded in $\sa$.





    Here we continue the study started in \cite{Dut06} of the support of the representations associated to super-wavelets. Starting with an embedding $\hat i$ of $\br^d$ into the solenoid, we have some periodic characters $\chi_i$ in $\sa$, associated to the cycles. Then we show that the representation is supported on the union of $\chi_i \hat i(\br^d)$ where the multiplication here is just the multiplication in $\sa$.

     We combine Mackey's theory of induced representations with the analysis of $\sa$-cycles. In this connection, we find a dynamical obstruction for embeddings as follows. Intuitively, one wants to encode the solenoid into a symbol space $[0,1)^d \times \Omega$. This can be a problem when one is dealing with matrices, as compared to the dyadic scaling in one dimension which is the traditional context for wavelet analysis.

The reason for the obstruction is that the candidate $[0,1)^d$ is not invariant under the inverse branches $\tau_d(x)=\at^{-1}(x+d)$, $d\in\mathcal D$ where $\mathcal D$ is a chosen finite set of vectors in $\br^d$. These maps serve as inverse branches to the action by $A^T$ on $\br^d/\bz^d$. The other candidate different from $[0,1)^d$ is the Hutchinson attractor $X(A^T,\mathcal D)$ for $\tau_d$, $X(A^T,\mathcal D)$ contained in $\br^d$. But $X(A^T,\mathcal D)$ might not tile $\br^d$ by $\bz^d$.

     So one question is: Can one choose $\mathcal D$ to be a complete set of representatives for the finite quotient group $\bz^d/ A^T\bz^d$ such that the attractor $X(A^T,\mathcal D)$ of $\{\tau_d\}$ tiles $\br^d$ by $\bz^d$?  If not, then how should one choose $A$ such that this is possible? Or perhaps, one must replace $A$ with an iterate $A^p$ of $A$ for a suitable $p$?
     
     The reader may find our use of the pair of scaling matrix $A$, and its transposed $A^T$ confusing. It is unavoidable and is dictated by our essential use of Fourier duality: If $A$ is acting in $d$-space, then $A^T$ is acting in the dual vector variable (say frequency), see e.g., Lemma \ref{lemfuglede}. In deciding tiling properties for $X = X(A^T,D)$ our use of spectral theory is essential as there seem to be no direct way of attacking the tiling problem for $X$.

    While initially Knuth's analysis \cite{Knu69} of what is now called affine Iterated Function Systems (IFS) was motivated by the desire to introduce geometry into algorithms for general digit sets in positional number systems, the idea of turning ``digits'' into geometry and tiling questions was followed up later by others, e.g., Odlyzko \cite{Odl78}, Hutchinson \cite{Hut81}, Bratteli-Jorgensen \cite{BrJo99}, and Lagarias-Wang \cite{LaWa96a, LaWa96b, LaWa96c, LaWa97, LaWa00}. In \cite{LaWa96a} the authors suggest that the tiling issues implied by the geometric positional ``number'' systems are directly connected with wavelet algorithms. In particular they pointed out that, for a fixed choice of $A$ and $\mathcal D$, the corresponding attractor $X(A^T,\mathcal D)$ as described above does not always tile $\br^d$ by translation vectors from the unit-grid lattice $\bz^d$. When it does, we say that $X(A^T,\mathcal D)$ is a Haar wavelet. The term ``Haar wavelet'' is used because the corresponding indicator function is a scaling function (father function) for an ONB wavelet system in $L^2(\br^d)$. The authors of \cite{LaWa96a} showed that in 2D, every expansive $d$ by $d$ matrix over $\bz$ has at least one ``digit'' set $\mathcal D$ such that $X(A^T,\mathcal D)$ is a Haar wavelet. It was later proved that in 5D, not every expansive 5 by 5 matrix over $\bz$ can be turned into a Haar wavelet; i.e., for such a matrix $A$, that there is no choice of $\mathcal D$ for which $X(A^T,\mathcal D)$ is a Haar wavelet.

      The fact that there are exotic 5 by 5 expansive matrices $A$ over $\bz$, i.e., $A$ in $\mathcal M_5(\bz)$ for which no digit set $\mathcal D$ may be found such that $X(A^T,\mathcal D)$ makes a $\bz^5$-tiling of $\br^5$ was worked out in the following papers \cite{LaWa96c, LaWa97, HLR02,HeLa04}. By digit set we mean a subset $\mathcal D$ in $\bz^5$, in bijective correspondence with $\bz^5/A^T\bz^5$. Such exotic matrices $A$ are said to not allow Haar wavelets.
    The question came up after Lagarias-Wang \cite{LaWa97} showed that every expansive $A$ in $\mathcal M_2(\bz)$ allows digit sets in $\bz^2$ which make $\bz^2$ tile, i.e, they allow Haar wavelets.

        The aim of this paper is to revisit the geometry of sets $X(A^T,\mathcal D)$ in light of recent results on IFS involving dynamics and representation theory, see e.g., \cite{Dut06, DuJo06a, DuJo06b, DuJo06c, DuJo07}.

\section{Definitions and notations}
While standard wavelet bases built on wavelet filters and on a fixed expansive $d$ by $d$ matrix over $\bz$, say $A$, refer to the Hilbert space $L^2(\br^d)$, many naturally occurring wavelet filters \cite{Dut06} suggest other Hilbert spaces than $L^2(\br^d)$, in fact Hilbert spaces containing a copy of $L^2(\br^d)$. This approach \cite{Dut06} suggests the the name ``super wavelets'', and naturally leads to representations of a Baumslag-Solitar group built on the matrix $A$, called wavelet representations. Starting with a fixed $A$, there is a compact solenoid $\sa$  with the property that matrix multiplication by $A$ on $\br^d/\bz^d$ induces an automorphism $\sigma_A$ on $\sa$.

\par In the past decade the literature on self-affine sets, encoding and
digit representations for radix matrices has grown; in part because of
applications to such areas as number theory, to dynamics, and to
combinatorial geometry. It is not possible here to give a complete list of
these directions. Our present work has been motivated by the papers \cite{AkSc05}, \cite{Cur06}, \cite{GaYu06}, \cite{HLR02}, \cite{HeLa04}, \cite{KLSW99}, \cite{Li06}, \cite{Li07}, \cite{Saf98},\cite{ZZ06}, \cite{LaWa00}.

\begin{definition}\label{deftaud}
Let $A$ be a $d\times d$ matrix with integer entries. We say that the matrix $A$ is {\it expansive} if all its eigenvalues $\lambda$ satisfy $|\lambda|>1$. 
\end{definition}

Let $A$ be a $d\times d$ expansive matrix with integer entries. Let
\begin{equation}\label{eqzda}
\bz^d[A^{-1}]=\{A^{-j}k\,|\,j\in\bn,k\in\bz^d\}.
\end{equation}

Note that $\bz^d[A^{-1}]$ is the inductive limit of the group inclusions
$$\bz^d\hookrightarrow A^{-1}\bz^d\hookrightarrow A^{-2}\bz^d\hookrightarrow\dots$$
or equivalently
$$\bz^d\stackrel{A}{\rightarrow}\bz^d\stackrel{A}{\rightarrow}\bz^d\stackrel{A}{\rightarrow}\dots$$

On $\bz^d[A^{-1}]$ we consider the discrete topology (even though $\bz^d[A^{-1}]$ is a subgroup of $\br^d$).

On the group $\bza$, the map $\alpha_A(x)=Ax$, $x\in\bza$ defines an automorphism of the group $\bza$. 

For the use of these groups in $C^*$-algebras, see e.g. \cite{BrJo99,BJKR01}.
\subsection{The group $G_A:=\bza\rtimes_{\alpha_A}\bz$}
The group $G_A:=\bza\rtimes_{\alpha_A}\bz$ is the semidirect product of $\bza$ under the action of $\bz$ by the automorphisms $\alpha_A$. This means that 
\begin{equation}\label{eqga}
G_A:=\{(j,b)\,|\,j\in \bz, b\in\bza\},\quad (j,b)\cdot(k,c)=(j+k,\alpha_A^j(c)+b),\quad(j,k\in\bz,b,c\in\bza).
\end{equation}
\begin{proposition}\label{propgagen}
The group $G_A$ is generated by the elements $u:=(1,0)$ and $t_k=(0,k)$, $k\in\bza$. Moreover
\begin{equation}
ut_ku^{-1}=t_{Ak},\quad(k\in\bz)
\end{equation}
\begin{equation}
t_{A^{-n}k}:=(0,A^{-n}k)=u^{-n}t_ku^n,\quad (j,A^{-n}k)=t_{A^{-n}k}u^j,\quad(n\geq 0,k\in\bz^d,j\in\bz).
\end{equation}
\end{proposition}

\begin{remark}
From Proposition \ref{propgagen} we infer that a unitary representation of the group $G_A$ is completely determined by giving some unitary operators $T_k$, $k\in\bz^d$ and $U$ subject to the relation 
\begin{equation}
UT_kU^{-1}=T_{Ak},\quad(k\in\bz).
\end{equation}
\end{remark}

 In \cite{DJ07b} we use induced representations of $G_A$ in the sense of
Mackey in order to encode wavelet sets for a fixed expansive $d$ by $d$
matrix over $\bz$. Note that Mackey's method was developed for continuous
groups, where encoding is done with co-adjoint orbits. In contrast we
show in the present paper that encoding in the solenoid is required
for wavelet representations, i.e., representations of discrete
versions of higher rank $ax + b$ groups.

\subsection{The dual group of $\bz^d[A^{-1}]$: the solenoid $\sa$.} The dual group of $\bz^d$ is $\bt^d$, where
$$\bt^d:=\{(z_1,\dots,z_d)\,|\, |z_i|=1, i\in\{1,\dots,d\}\}.$$
For $x=(x_1,\dots,x_d)\in\br^d$ let 
$$e^{2\pi ix}:=(e^{2\pi ix_1},\dots,e^{2\pi ix_d})\in\bt^d.$$
For $k=(k_1,\dots,k_d)\in\bz^d$ and $z=(e^{2\pi ix_1},\dots,e^{2\pi ix_d})\in\bt^d$, we use the notation
$$z^k:=e^{2\pi ik_1x_1+\dots +2\pi ik_dx_d}=e^{2\pi ik\cdot x}\in\bt.$$

The duality between $\bz^d$ and $\bt^d$ is given by
\begin{equation}\label{eqdualzt}
\ip{k}{z}=z^k,\quad(k\in\bz^d,z\in\bt^d).
\end{equation}

For $z=(e^{2\pi ix_1},\dots e^{2\pi ix_d})\in\bt^d$ we denote by 
\begin{equation}\label{eqza}
z^{A}:=e^{2\pi i\at x}=(e^{2\pi i \sum_{j=1}^da_{j1}x_j},\dots,e^{2\pi i\sum_{j=1}^da_{jd}x_d}).
\end{equation}
Note that $z^{Ak}=(z^A)^k$ for all $k\in\bz^d$. 
\begin{definition}\label{defsole}
The dual group of $\bza$ is the group $\sa$ defined by
$$\sa:=\{(z_n)_{n\in\bn}\,|\, z_n\in\bt^d, z_{n+1}^A=z_n,\mbox{ for all }n\in\bz\}.$$
The group $\sa$ is called the {\it solenoid} of $A$. It is a compact Abelian group with the topology induced from the product topology on $\bt^\bn$.  
\end{definition}

The duality is given by
\begin{equation}\label{eqdualza}
\ip{A^{-j}k}{(z_n)_{n\in\bn}}=\ip{k}{z_j}=z_j^k,\quad(j\in\bn,k\in\bz^d,(z_n)_{n\in\bn}\in\sa).
\end{equation}

 The dual of the automorphism $\alpha_A$ on $\bza$, $\alpha_A(x)=Ax$ is the {\it shift}
\begin{equation}\label{eqshift}
\sigma_A(z_0,z_1,\dots)=(z_0^A,z_0,z_1,\dots),\quad((z_0,z_1,\dots)\in\sa).
\end{equation}

We denote by $\theta_n$ the projection maps $\theta_n:\sa\rightarrow\bt^d$, $\theta_n(z_0,z_1,\dots)=z_n$ for $n\in\bn$. Note that
\begin{equation}\label{eqtetan}
\theta_{n+1}\circ\sigma_A=\theta_n,\quad \left(\theta_{n+1}(z_0,z_1,\dots)\right)^A=\theta_n(z_0,z_1,\dots),\quad(n\in\bn,(z_0,z_1,\dots)).
\end{equation}

\section{Embeddings of $\br^d$ into the solenoid $\sa$}\label{embed}
\par We saw in Proposition \ref{propgagen} that there is a natural semidirect product discrete group $G_A$ which carries a unitary wavelet representation. For wavelets in $\br^d$, this construction begins with a fixed expansive $d$ by $d$  matrix $A$, and the unitary representation will be acting in the Hilbert space $L^2(\br^d)$. It is known that a certain redundancy (\cite{BDP05}, \cite{HaLa00}) in wavelet constructions dictates unitary representations in Hilbert spaces larger than $L^2(\br^d)$, i.e., with $L^2(\br^d)$ embedded as an isomorphic copy in an ambient Hilbert space. By introducing a specific embedding of $\br^d$ in an ambient solenoid $S_A$ we are able to account for super representations (Definition \ref{c1}).  The action of $A$ induces an automorphism $\sigma_A$ in $S_A$. By computing periodic points for  $\sigma_A$  we are able (Theorem \ref{prop3_9}) to account for the super representations acting in an $L^2$ space defined from an induced measure on $S_A$. 
\par In sections 5 and 6 we will further study periodic points and cycles. It turns out that the notion of cycle is different when referring to the integer points and the fractions. Starting with a fixed radix pair $(A,\mathcal D)$, we will make ``integer points'' precise in terms of associated lattices (rank-d subgroups in $\br^d$), and the ``fractions'' will take the form of compact subsets $X$ in $\br^d$ defined by an $(A,\mathcal D)$ self-similarity (Definition \ref{ifs}).

 We begin by showing how the space $\br^d$ can be seen as subspace of the solenoid $\sa$.
\begin{proposition}\label{propembed1}
The inclusion $i:\bza\rightarrow\br^d$ has a dual $\hat i:\br^d\rightarrow\sa$
\begin{equation}\label{eqhati}
\hat i(x)=(e^{2\pi i\at^{-n}x})_{n\in\bn},\quad(x\in\br^d).
\end{equation}

The map $\hat i$ is one-to-one, and onto the set of sequences $(z_n)_{n\in\bn}\in\sa$ with the property that $\lim_{n\rightarrow\infty}z_n=\mathbf 1$ where $\mathbf 1=(1,\dots,1)$ is the neutral element of $\bt^d$.

The map $\hat i$ satisfies the following relation
\begin{equation}\label{eqhatia}
\hat i(A^Tx)=\sigma_A(\hat i(x)),\quad(x\in\br^d).
\end{equation}
\end{proposition}

\begin{proof}
To see that $\hat i$ is one-to-one, we notice that if $x,x'\in\br^d$ and $\hat i(x)=\hat i(x')$ then 
$\at^{-n}x-\at^{-n}x'\in\bz^d$ for all $n\in\bn$. Since $A$ is expansive, the norm $\|\at^{-n}(x-x')\|$ converges to $0$ as $n\rightarrow\infty$. Thus, $x$ and $x'$ must be equal.

Since $A$ is expansive, so is $A^T$, so $\at^{-n}x$ converges to $0$ for all $x\in\br^d$. Therefore $e^{2\pi iA^{-n}x}$ converges to $\mathbf 1$. Conversely, suppose $(z_n)_{n\in\bn}$ is in $\sa$ and $z_n$ converges to $1$. Then
$z_n=e^{2\pi ix_n}$ for some $x_n\in\br^d$ and, for $n$ large we can assume $x_n$ is close to $0$. Since we have $z_{n+1}^A=z_n$, this implies that $A^Tx_{n+1}\equiv x_n\mod\bz^d$, so $A^Tx_{n+1}=x_n+l$ for some $l\in\bz^d$. But since both $x_{n+1}$ and $x_n$ are close to $0$, this implies that, for some $n_0\in\bn$, we have $\at x_{n+1}=x_n$ for $n\geq n_0$. Let $x:=\at^{n_0} x_{n_0}$. The previous argument shows that $\at^jx_{n_0}\equiv x_{n_0-j}\mod\bz^d$ for all $j\leq n_0$, so $\hat i(x)=(e^{2\pi i\at^{-n}x})_{n\in\bn}$. Thus, $\hat i$ is onto.

The other assertions follow from some direct computations, using the duality in \eqref{eqdualza}.
\end{proof}
 Now that we have the embedding of $\br^d$ into the solenoid, we can transport the wavelet representation on $L^2(\br^d)$ to the solenoid $\sa$.
\begin{definition}\label{deftul2}
On $L^2(\br^d)$ we denote by $T_k$ the {\it translation operator} $(T_kf)(x)=f(x-k)$, $g\in L^2(\br^d)$, $x\in\br^d$, $k\in\bz^d$, and by $U$ the {\it dilation operator} $(Uf)(x)=\frac{1}{\sqrt{|\det A|}}f(A^{-1}x)$. Their Fourier transform is
\begin{equation}\label{eqhattu}
(\hat T_kh)(x)=e^{2\pi ik\cdot x}f(x),\quad(\hat Uh)(x)=\sqrt{|\det A|}h(A^Tx),\quad(h\in L^2(\br^d),x\in\br^d,k\in\bz^d).
\end{equation}
The operators $\{U,T_k\}$ (or $\{\hat U,\hat T_k\}$) define a representation of the group $G_A$ on $L^2(\br^d)$.
\end{definition}
\begin{definition}
We denote by $\sa(1)$ the set of sequences $(z_n)_{n\in\bn}\in\sa$ such that $\lim_{n\rightarrow\infty}z_n=\mathbf 1$. 
On $\sa(1)$ consider the measure $\tilde\mu$ defined by
\begin{equation}\label{eqmutilda}
\int_{\sa(1)}f\,d\tilde\mu=\int_{\bt^d}\sum_{(z_n)_n\in\sa(1),\theta_0((z_n)_{n\in\bn})=z}f((z_n)_{n\in\bn})\,d\mu(z),
\end{equation}
Where $\mu$ is the Haar measure on $\bt^d$. 

On $L^2(\sa(1),\tilde\mu)$ define the operators 
\begin{equation}\label{eqTtilda}
(\tilde T_kf)(z_0,z_1,\dots)=z_0^kf(z_0,z_1,\dots),\quad(z_0,z_1,\dots)\in\sa(1),k\in\bz^d),
\end{equation}
\begin{equation}\label{eqUtilda}
(\tilde Uf)(z_0,z_1,\dots)=\sqrt{|\det A|}f(\sigma_A(z_0,z_1,\dots)),\quad((z_0,z_1,\dots)\in\sa(1)).
\end{equation}
\end{definition}

   Our next theorem shows that this unitary resentation of the
reduced solenoid   $\sa(1)$ is a universal super representation of the
wavelet group in that when $A$ is given, then via an intertwining
isometry $\mathcal W: L^2(\br^d)\rightarrow L^2(\sa(1))$ the standard $A$-wavelet-unitary
representation acting on $L^2(\br^d)$ is naturally included in the $\sa(1)$
representation. The details of the symbolic encoding of this pair of representations depends on
a choice of digit set $\mathcal D$, and the structure of the associated
$(A,\mathcal D)$-cycles (we will give these details in Section 6).

\begin{theorem}
(i) The measure $\tilde\mu$ satisfies the following invariance property
\begin{equation}\label{eqinvmutilda}
\int_{\sa(1)}f\circ\sigma_A\,d\tilde\mu=\frac1{|\det A|}\int_{\sa(1)}f\,d\tilde\mu,\quad(f\in L^1(\sa(1),\tilde\mu)).
\end{equation}

(ii) The operators $\tilde T_k$, $k\in\bz^d$ and $\tilde U$ are unitary and they satisfy the following relation
\begin{equation}\label{eqcov}
\tilde U\tilde T_k\tilde U^{-1}=\tilde T_{Ak},\quad(k\in\bz^d)
\end{equation}
so $\{\tilde U,\tilde T_k\}$ generate a representation of the group $G_A$ on $L^2(\sa(1),\tilde\mu)$.

(iii) The map $\hat i$ is a measure preserving transformation between $\br^d$ and $\sa(1)$.

(iv) The operator $\mathcal W:L^2(\br^d)\rightarrow L^2(\sa(1),\tilde\mu)$ defined by $\mathcal Wf=f\circ\hat i^{-1}$ is an intertwining isometric isomorphism,
\begin{equation}\label{eqhatiinter}
\mathcal W\hat T_k=\tilde T_k\mathcal W,\quad(k\in\bz^d),\quad \mathcal W\hat U=\tilde U\mathcal W.
\end{equation}

\end{theorem}

\begin{proof}
(i) The Haar measure on $\bt^d$ satisfies the following strong invariance property:
\begin{equation}\label{eqinvmu}
\int_{\bt^d}f\,d\mu=\frac{1}{|\det A|}\int_{\bt^d}\sum_{y^A=z}f(y)\,d\mu(z),\quad(f\in L^1(\mu)).
\end{equation}
Using this, we have for $f\in L^1(\sa(1),\tilde\mu)$:
$$\int_{\sa(1)}f\circ\sigma_A\,d\tilde\mu=\int_{\bt^d}\sum_{(z_n)_n\in\sa(1),z_0=z}f(z_0^A,z_0,z_1,\dots)\,d\mu(z)=$$
$$\frac1{|\det A|}\int_{\bt^d}\sum_{y^A=z}\sum_{(z_n)_n\in\sa(1), z_0=y}f(z_0^A,z_0,z_1,\dots)\,d\mu(z)=
\frac{1}{|\det A|}\int_{\bt^d}\sum_{(w_n)_n\in\sa(1),w_0=z}f(w_0,w_1,\dots)\,d\mu(z)=$$$$\frac{1}{|\det A|}\int_{\sa(1)}f\,d\tilde\mu.$$

(ii) follows from (i) and some direct computations.

(iii) First note that, by Proposition \ref{propembed1} $(z_n)_n\in\sa(1)$ with $z_0=e^{2\pi ix}$ iff $(z_n)_n=\hat i(y)$ and $e^{2\pi iy}=e^{2\pi ix}$, i.e., $(z_n)_n=\hat i(x+k)$ for some $k\in\bz^d$. 

Let $f\in L^1(\sa(1),\tilde\mu)$. Then 
$$\int_{\br^d}f(\hat i(x))\,dx=\int_{[0,1)^d}\sum_{k\in\bz^d}f(\hat i(x+k))\,dx=$$
$$\int_{[0,1)^d}\sum_{(z_n)_n\in\sa(1),z_0=e^{2\pi ix}}f(z_0,z_1,\dots)\,dx=\int_{\bt^d}\sum_{(z_n)_n\in\sa(1),z_0=z}f(z_0,z_1,\dots)\,d\mu(x)=\int_{\sa(1)}f\,d\tilde\mu.$$
This proves that $\hat i$ is measure preserving.

(iv) Since $\hat i$ is measure preserving, $\mathcal W$ is an isometric isomorphism. The intertwining relation \eqref{eqhatiinter} follows by a direct computation that uses \eqref{eqhatia}.
\end{proof}

\subsection{Cycles}
 Next we will show how the ``super-wavelet'' representations from \cite{BDP05} can be realized on the solenoid.
\begin{definition}\label{c1}
An ordered set $C:=\{\zeta_0,\zeta_1,\dots,\zeta_{p-1}\}$ in $\bt^d$ is called a {\it cycle} if $\zeta_{j+1}^A=\zeta_j$ for $j\in\{0,\dots,p-2\}$ and $\zeta_0^A=\zeta_{p-1}$, where $p\geq 1$. The number $p$ is called the {\it period} of the cycle if the points $\zeta_i$ are distinct.  
\end{definition}

\begin{definition}\label{defhc}
Let $C=\{\zeta_0,\zeta_1,\dots,\zeta_{p-1}\}$ be a cycle. Denote by 
\begin{equation}\label{eql2rlap}
\mathcal H_C:=\underbrace{L^2(\br^d)\oplus\dots\oplus L^2(\br^d)}_{p\mbox{ times}}=L^2(\br^d\times\bz_p),
\end{equation}
where $\bz_p=\{0,\dots,p-1\}$ is the cyclic group of order $p$.

Define the operators on $\mathcal H_C$:
\begin{equation}\label{eqtc}
T_{C,k}(f_0,\dots,f_{p-1})=(\zeta_0^kT_kf_0,\dots,\zeta_{p-1}^kT_kf_{p-1}),\quad((f_0,\dots,f_{p-1})\in\mathcal H_C,k\in\bz^d)
\end{equation}

\begin{equation}\label{equc}
U_C(f_0,\dots,f_{p-1})=(Uf_{p-1},Uf_0,\dots,Uf_{p-2}),\quad((f_0,\dots,f_{p-1})\in\mathcal H_C)
\end{equation}
where $T_k$ and $U$ are the operators on $L^2(\br^d)$ from Definition \ref{deftul2}. We will denote by $\hat T_{C,k}$ and $\hat U_C$ the Fourier transform of these operators (the Fourier transform being applied on each component of $\mathcal H_C$). 
\end{definition}

Then a simple calculations shows that:
\begin{proposition}
The operators $T_{C,k}$, $k\in\bz^d$ and $U_C$ are unitary and satisfy the following relation
\begin{equation}\label{eqcovc}
U_CT_{C,k}U_C^{-1}=T_{C,Ak},\quad(k\in\bz^d)
\end{equation}
so they define a representation of the group $G_A$ on $\mathcal H_C$.
\end{proposition}

Some examples and wavelet applications are also included in \cite{Jor03}.

\begin{definition}
Let $C:=\{\zeta_0,\dots,\zeta_{p-1}\}$ be a cycle, $\zeta_j=e^{2\pi i\theta_j}$, for some $x_j\in\br^d$, $(j\in\{0,\dots,p-1\})$. We will use the notation $\theta_n:=\theta_{n\mod p},\zeta_n:=\zeta_{n\mod p}$ for all $n\in\bz$. We denote by 
$$\chi_C:=(\zeta_0,\dots,\zeta_{p-1},\zeta_0,\dots,\zeta_{p-1},\dots)\in\sa.$$
Let 
\begin{equation}
	\sa(C):=\bigcup_{j=0}^{p-1}\sigma_A^{-j}(\chi_C)\hat i(\br^d).
	\label{eq:}
\end{equation}
Let $\hat i_C:\br^d\times\bz_p\rightarrow\sa(C)$
\begin{equation}
	\hat i_C(x,j)=\sigma_A^{-j}(\chi_C)\hat i(x)=(e^{2\pi i(\at^{-n}x+\theta_{n+j})})_{n\in\bn},\quad(x\in\br^d,j\in\bz_p).
	\label{eqic}
\end{equation}
Define the measure $\tilde\mu$ on $\sa(C)$ by an equation similar to \eqref{eqmutilda}, (the only difference here is the support $\sa(C)$ instead of $\sa(1)$)
$$\int_{\sa(C)}f\,d\tilde\mu=\int_{\bt^d}\sum_{(z_n)_n\in\sa(C),\theta_0((z_n)_{n\in\bn})=z}f((z_n)_{n\in\bn})\,d\mu(z).
$$
Define the operators $\tilde T_k$, $k\in\bz^d$ and $\tilde U$ on $L^2(\sa(C),\tilde\mu)$ by the same formulas as in \eqref{eqTtilda} and \eqref{eqUtilda}.
\end{definition}

\begin{theorem}\label{prop3_9}
(i) The point $\chi_C$ is periodic for $\sigma_A$, $\sigma_A^p(\chi_C)=\chi_C$, and $\sigma_A$ permutes cyclically the sets $\sigma^{-j}(\chi_C)\hat i(\br^d)$, $j\in\bz_p$. 

(ii) The set $\sa(C)$ consists of exactly the points $(z_n)_{n\in\bn}$ with the property that the distance from $z_n$ to $C$ converges to $0$.

(iii) Let $\alpha_{A,p}:\bz^d\times\bz_p\rightarrow \bz^d\times\bz_p$
\begin{equation}\label{eqalphap}
\alpha_{A,p}(x,j)=(A^Tx,j-1),\quad(x\in\br^d,j\in\bz_p).
\end{equation}
The map $\hat i_C$ is a bijective measure preserving transformation that satifies the relation
\begin{equation}
	\hat i_C\circ\alpha_{A,p}=\sigma_A\circ\hat i_C.
	\label{eqica}
\end{equation}

(iv) The operator $\mathcal W_C: \mathcal H_C=L^2(\br^d\times\bz_p)\rightarrow L^2(\sa(C),\tilde\mu)$, $\mathcal W_Cf=f\circ\hat i_C^{-1}$, is an intertwining isometric isomorphism:
\begin{equation}\label{eqwcint}
\mathcal W_C\hat T_{C,k}=\tilde T_k\mathcal W_C,\quad(k\in\bz^d),\quad\mathcal W_C\hat U_C=\tilde U\mathcal W_C.
\end{equation}
\end{theorem}

\begin{proof}
(i) is trivial. 

(ii) If $(z_n)_n$ is in $\sa(C)$ then for some $j$, we have that $(\sigma_A^{-j}(\chi_C))^{-1}(z_n)_n$ is in $\hat i(\br^d)$ so 
$\zeta_j^{-1}z_{np}$ converges to $0$. Therefore $z_{np}$ converges to $\zeta_j$ so $z_{np-l}=z_{np}^{A^l}$ converges to $\zeta_{j-l}$ for all $l\in\bz_p$.

Conversely, suppose $(z_n)_n\in\sa$ and $\dist(z_n,C)$ converges to $0$. We claim that $z_{np}$ converges to one of the points of the cycle $\zeta_i$. 

Pick an $\epsilon>0$ small enough such that 
for all $i\in\bz_p$, $\dist(z_n,\zeta_i)<\epsilon$ implies $\dist(z_n,\zeta_{i'})>\epsilon$ for $i\neq i'$. There exists a $\delta>0$, $\delta<\epsilon$ such that for all $i\in\bz_p$, if $\dist(z,\zeta_i)<\delta$ then $\dist(z^{A^p},\zeta_i)<\epsilon$. 

There exists an $n_\epsilon$ such that if $n\geq n_\epsilon$ then $\dist(z_n,C)<\delta<\epsilon$. 
Then for some $i\in\bz_p$ we have $\dist(z_{n_\epsilon},\zeta_i)<\epsilon$. Also for some $i'\in\bz_p$ we have $\dist(z_{n_\epsilon+p},\zeta_{i'})<\delta$. This implies that $\dist(z_{n_\epsilon},\zeta_{i'})<\epsilon$ so $i'=i$. 
By induction, we obtain that $\dist(z_{n_\epsilon+kp},\zeta_i)<\delta<\epsilon$. And since $\dist(z_n,C)$ converges to $0$ this shows that $\dist(z_{n_\epsilon+kp},\zeta_i)$ converges to $0$. Applying the map $z^A$ several times we conclude that $z_{np}$ converges to one of the elemenst of the cycle, $\zeta_j$. 

Then consider $(w_n)_n:=\sigma_A^j(\chi_C)^{-1}(z_n)_n\in\sa$. Clearly $w_n$ converges to $\mathbf 1$. By Proposition \ref{propembed1}, there exists an $x\in\br^d$ such that $(w_n)_n=\hat i(x)$. Thus $(z_n)_n=\sigma_A^j(\chi_C)\hat i(x)$, and this proves (ii).

(iii) Since $\hat i$ is bijective (Proposition \ref{propembed1}), clearly $\hat i_C$ is also bijective. To check that $\hat i_C$ is a measure preserving transformation, since $\hat i$ is measure preserving by Proposition \ref{propembed1}, it is enough to check that multiplication by $\sigma_A^{-j}(\chi_C)$ leaves the measure $\tilde\mu$ invariant, i.e., for a function $f$ defined on $\sigma_A^{-j}(\br^d)$,
\begin{equation}\label{eqinvtrmutilda}
\int_{\sigma_A^{-j}(\chi_C)\hat i(\br^d)}f\,d\tilde\mu=\int_{\hat i(\br^d)}f(\sigma_A^{-j}(\chi_C)(z_n)_n)\,d\tilde\mu(z_n)_n,\quad(j\in\bz_p).
\end{equation}
It is enough to check this for $j=0$. Using the translation invariance of the Haar measure $\mu$ on $\bt^d$, we have:
$$\int_{\chi_C\hat i(\br^d)}f\,d\tilde\mu=\int_{\bt^d}\sum_{(z_n)_n\in\chi_C\hat i(\br^d),z_0=z}f((z_n)_n)\,d\mu(z)=
\int_{\bt^d}\sum_{(w_n)_n\in\hat i(\br^d),w_0=z\zeta_0^{-1}}f(\chi_C(w_n)_n)\,d\mu(z)=$$
$$\int_{\bt^d}\sum_{(w_n)_n\in\hat i(\br^d),w_0=z}f(\chi_C(w_n)_n)\,d\mu=\int_{\hat i(\br^d)}f(\chi_C(w_n)_n)\,d\tilde\mu(w_n)_n.$$

Equation \eqref{eqica} follows by a direct computation.

(iv) Follows from \eqref{eqica}.
\end{proof}

\section{Encoding of integer points}

       The idea of using matrices and geometry in creating a positional number system for points in $\br^d$ was initiated by Don Knuth, see especially \cite{Knu76}, vol 2, chapter 4 (Arithmetic), and 4.1 which introduces this geometric and algorithmic approach to positional number system. In fact, the Twin-Dragon appears on page 206 (in v 2 of \cite{Knu76}).  In our present discussion, with a fixed expansive matrix $A^T$ playing the role of the basis-number (or the radix) in our radix representations, the natural question arises: ``What is the role of the integer lattice $\bz^d$ relative to our radix system?'' This section gives a preliminary answer to the question, and the next section is a complete analysis involving cycles.

       As before, we begin with a choice of expansive matrix $A$ (i.e., $A$ is a fixed $d$ by $d$ matrix over $\bz$), and we choose a subset $\mathcal D$ in $\bz^d$ (for digits), points in $\mathcal D$ in bijective correspondence with $\bz^d/A^T\bz^d$. As it turns out, ``the integers'' relative to the $(A,\mathcal D)$-radix typically will not have finite radix (or Laurent) expansions in positive powers of $A^T$. The reason for this is the presence of certain non-trivial cycles $C$ in $\bz^d$ leading to infinite repetitions, which we will take up systematically in the next section.

      In this section we will show that when the pair $(A,\mathcal D)$ is fixed, there is an encoding mapping $\phi$ which records the finite words in the alphabet $\mathcal D$ which will correspond to the cycles in $\bz^d$ that are associated with our particular choice of $(A,\mathcal D)$. However, once the cycles in $\bz^d$ are identified, there is a much more detailed encoding directly for $\bz^d$ which will be done in detail in Theorem \ref{thenccycl}.

      So our present encoding mapping  $\phi : \bz^d\rightarrow  \mathcal D^{\bn}$, depending on the pair $(A,\mathcal D)$, is an introduction to our analysis of a refined solenoid encoding and of all cycles in the next section. Still our starting point is a fixed radix-pair $(A,\mathcal D)$. The fact that the encoding mapping  $\phi : \bz^d\rightarrow \mathcal D^{\bn}$ is injective is a consequence of the expansive property. Corresponding to $(A,\mathcal D)$ there is a finite set of finite words $F = F(A,\mathcal D)$ in letters from $\mathcal D$. These words label the $\bz^d$-cycles, and the encoding mapping  $\phi : \bz^d\rightarrow \mathcal D^{\bn}$ (infinite Cartesian product) maps onto the set of infinite words which terminate in an infinite repetition of one of the words from $F$.

Let $d$ be given. Let $B$ be a $d\times d$ matrix over $\bz$, and assume 
\begin{equation}\label{eqdz1}
\bigcap_{k\geq 1} B^k\bz^d=\{0\}.
\end{equation}
Note that this holds if $B$ is assumed expansive.

Let $\mathcal D\subset\bz^d$ be a complete set of representatives for $\bz^d/B\bz^d$. Assume $0\in\mathcal D$. (This assumption is for convenience and can be easily removed {\it mutatis mutandis}.)
\begin{definition}
Define a $(B,\mathcal D)$ encoding of $\bz^d$  
$$\phi(x):=d_0d_1d_2\dots, \phi:\bz^d\rightarrow\mathcal D^{\bn}$$
as follows: when $x\in\bz^d$ is given, there is a unique pair $x_1\in\bz^d$ and $d_0\in\mathcal D$ such that
\begin{equation}\label{eqdz3}
x=d_0+Bx_1.
\end{equation}
By the same argument, now determine $d_1,d_2,\dots\in\mathcal D$, and $x_2,x_3,\dots\in\bz^d$ recursively such that
\begin{equation}\label{eqdz4}
x_k=d_k+Bx_{k+1}.
\end{equation}
\end{definition}
\begin{definition}
Set $$\Omega:=\mathcal D^{\bn}=\prod_{n=0}^\infty \mathcal D.$$
Elements $\omega\in\Omega$ are called infinite words in the alphabet $\mathcal D$. If $v$ is a finite word, we denote by $\underline v$ the infinite repetition of this word $vvv\dots$. If there are finite words $v$ and $w$ such that $\omega=(v\underline w)$, we say that $\omega$ ends in a cycle. 
\end{definition}
\begin{proposition}\label{propdzdz}
(i) The encoding mapping $\phi:\bz^d\rightarrow\Omega$  is well defined.

(ii) $\phi$ is one-to-one.

(iii) $\phi$ maps onto a subset of $\Omega$ of all infinite words that end in cycles. 
\end{proposition}

\begin{proof}
Part (i) is immediate from \eqref{eqdz3}, \eqref{eqdz4}.

(ii) Suppose $x,y\in\bz^d$ and $\phi(x)=\phi(y)$. Then an application of \eqref{eqdz4} shows that $x-y\in\cap_{k\geq 1}B^k\bz^d$, and we conclude that $x=y$ by an application of \eqref{eqdz1}.

(iii) Follows immediately from Theorem \ref{thenccycl}.
\end{proof}

\begin{remark}
(i) For $x\in\bz^d$, the encoding $\phi(x)=v\underline w=d_0d_1d_2\dots$ with $d_i\in\mathcal D$, and $v,w$ finite  words is unique; but the formal sum 
\begin{equation}\label{Eqdz5}
d_0+Bd_1+B^2d_2+\dots
\end{equation}
is not convergent unless $\underline w=\underline 0=000\dots$ infinite repetition. In that case there exists $m\in\bn$ such that $v=d_0\dots d_{m-1}$ and 
$$x=d_0+Bd_1+\dots+B^{m-1}d_{m-1}.$$

(ii) Suppose $v=\emptyset$ and $\phi(x)=\underline w$, with $w=l_0l_1\dots l_{p-1}$. Then $-x$ has the following infinite, convergent, periodic fractional expansion 
$$-x=\sum_{k=0}^\infty\sum_{i=0}^{p-1} B^{-kp+i+1}l_{p-1-i}.$$
See Proposition \ref{proprc} and Theorem \ref{thenccycl} for the proof.
\end{remark}

\begin{example}
The following simple example in 1D illustrates the cases (i) and (ii) above.

Let $d=1$ , $B=2$, and $\mathcal D=\{0,3\}$. Let $\phi:\bz\rightarrow\Omega$ be the encoding. 

We have $\phi(11)=3003\underline{30}$, i.e., $v=3003$ and $w=30$. 

And $\phi(18)=033\underline0$, i.e., $v=033$ and $w=0$, corresponding to the finite representation 
$$18=0+3\cdot 2+3\cdot 2^2.$$
Finally $\phi(-2)=\underline{03}$, i.e., $v=\emptyset$, and $w=03$. Hence by (ii) we get the following infinite fractional dyadic representation $$2=3\cdot 2^{-1}+0\cdot 2^{-2}+3\cdot 2^{-3}+0\cdot 2^{-4}+3\cdot 2^{-5}+\dots$$

Proposition \ref{proprc} shows that the cycles are obtained by intersecting $\bz$ with the set $-X(B,\mathcal D)$, where $X(B,\mathcal D)$ is the attractor of the maps $\tau_0(x)=x/2$, $\tau_3(x)=(x+3)/2$. In our example $X(B,\mathcal D)=[0,3]$. There are two cycles of length one $\{0\}$ and $\{-3\}$ and one cycle of length two $\{-1,-2\}$.

The encoding mapping $\phi$ records the cycles as follows: the one-cycles $\phi(0)=\underline 0$, $\phi(-3)=\underline 3$. The two-cycle: $\phi(-1)=\underline{30}$, $\phi(-2)=\underline{03}$. 
\end{example}
The next section takes up the encodings in general.

\section{Encodings of the solenoid}
 In this section we return to the geometry of general digit sets in positional number systems, turning ``digits'' into geometry and tilings. The starting point is a given pair $(A,\mathcal D)$ with $A$ expansive over $\bz$, and $\mathcal D$ a complete digit set. With the aid of the solenoid we give an explicit encoding. Specifically, we show that the attractor $X(A^T,\mathcal D)$ for the corresponding affine Iterated Function System (IFS) is a set of fractions for an $(A,\mathcal D)$-digital representation of points in $\br^d$. Moreover our positional ``number representation'' is spelled out in the form of an explicit IFS encoding of the compact solenoid $\sa$ associated with the pair $(A,\mathcal D)$.  The intricate part (Theorem \ref{thenccycl}) is played by the cycles in $\bz^d$ for the initial $(A,\mathcal D)$-IFS. Using the cycles we are able to write down formulas for the two maps which do the encoding as well as the decoding in our positional $\mathcal D$-representation.

Take a point $(z_n)_{n\in\bn}$ in the solenoid $\sa$. Then $z_0\in\bt^d$, $z_1^A=z_0$, $z_2^A=z_1$, and so on.

Since $z_0$ is in $\bt^d$ it can represented by $e^{2\pi ix_0}$, where $x_0\in\br^d$. Note that one has several choices for $x_0$, any of its integer translates $x_0+k$, $k\in\bz^d$, will do. 

Then $z_1^A=z_0$, so $z_1$ is a ``root'' of $z_0$. There are $|\det A|$ choices: if $z_0=e^{2\pi ix_0}$ then $z_1$ must be one of the points $e^{2\pi i\at^{-1}(x_0+d)}$, $d\in\mathcal D$, where $\mathcal D$ is some complete set of representatives for $\bz^d/A^T\bz^d$. Say $z_1=e^{2\pi i\at^{-1}(x_0+d_0)}=e^{2\pi ix_1}$.

At the next step $z_2$ is a root of $z_1$ so $z_2=e^{2\pi i \at^{-1}(x_1+d_1)}$ for some $d_1\in\mathcal D$. By induction, we get a sequence $d_1,d_2,\dots$ in $\mathcal D$. 

Thus, picking a point $(z_n)_{n\in\bn}$ in $\sa$ amounts to choosing an $x_0\in\br^d$ and an infinite word $d_0d_1\dots\in\Omega:=\mathcal D^{\bn}$. Thus we say that $(z_n)_{n\in\bn}$ can be {\it encoded} as 
$$(z_n)_{n\in\bn}\leftrightarrow (x_0,d_0d_1\dots).$$

Now note that changing the choice of $x_0$ (to say $x_0+k$), affects the entire sequence $d_0,d_1,\dots$. We want to make this choice unique in some sense. For this we need to find a subset $F$ of $\br^d$, such that for each $z\in\bt^d$, there is a unique $x\in F$ such that $z=e^{2\pi ix}$. In other words, $F$ must tile $\br^d$ by $\bz^d$-translations. 

Of course a first choice of this set $F$ would be $[0,1)^d$. While this works in dimension $d=1$, it may be inappropriate for higher dimensions. The problem is that we would need also $z_1$ to come from $e^{2\pi ix_1}$ with $x_1\in F$, and this would mean that $x_1=\at^{-1}(x_0+d_0)$ is in $F$. Thus, our set $F$ must have the following property 
$$\bigcup_{d\in\mathcal D}\at^{-1}(F+d)\subset F.$$
But since $(z_1,z_2,\dots)$ is also an element of $\sa$ and $z_1$ can be any point in $\bt^d$, it follows that we must actually have
\begin{equation}\label{eqF}
\bigcup_{d\in\mathcal D}\at^{-1}(F+d)= F.
\end{equation}
Of course, when we are interested only in measure theoretic notions, we can allow the equalities to hold only almost everywhere.

When $F$ is compact, this equation identifies $F$ as the attractor of an affine iterated function system.
\begin{definition}\label{ifs}
Let $\mathcal D$ be a complete set of representatives for $\bz^d/A^T\bz^d$. For every $d\in\mathcal D$, we denote by $\tau_d$ the map on $\br^d$ defined by 
\begin{equation}\label{eqtaud}
\tau_d(x)=\at^{-1}(x+d),\quad(x\in\br^d).
\end{equation}
\end{definition}

With this notation, if $F$ is compact, equation \eqref{eqF} says that $F$ is the attractor of the affine iterated function system $(\tau_d)_{d\in\mathcal D}$. This identifies $F$ as 
\begin{equation}\label{eqxa}
F=X(A^T,\mathcal D):=\{\sum_{j=1}^\infty \at^{-j}d_j\,|\,d_j\in\mathcal{D}\}.
\end{equation}

To find an encoding of the solenoid $\sa$ means to find a subset $F$ of $\br^d$ and a complete set of representatives $\mathcal D$ of $\bz^d/A^T\bz^d$ such that, if $x\in F$ then $\tau_dx\in F$ for all $d\in\mathcal D$, and the {\it decoding} map $\mathfrak d:F\times\mathcal{D}^\bn\rightarrow \sa$ defined by
$$F\times\mathcal D^{\bn}\ni(x_0,d_0d_1\dots)\mapsto (e^{2\pi ix_0},e^{2\pi i\tau_{d_0}x_0},e^{2\pi i\tau_{d_1}\tau_{d_0}x_0},\dots)\in\sa$$
is a bijection. 

Thus, to find this encoding of $\sa$ we need a subset $F$ of $\br^d$ and a complete set of representatives $\mathcal D$ that satisfy \eqref{eqF} and such that $F$ tiles $\br^d$ by integer translations, i.e., 
\begin{equation}\label{eqtile}
\bigcup_{k\in\bz^d}(F+k)=\br^d\mbox{ Lebesgue-a.e., and }(F+k)\cap(F+k')=\emptyset\mbox{ Lebesgue-a.e.}
\end{equation}
 
So the problem of encoding the solenoid into the space $F\times\mathcal D^{\bn}$ is equivalent to the following:

{\bf Question.} Given an expansive $d\times d$ integer matrix $A$, is it possible to find a complete set of representatives $\mathcal D$ of $\bz^d/A^T\bz^d$ such that the attractor $X(A^T,\mathcal D)$ of the iterated function system $(\tau_d)_{d\in\mathcal D}$ tiles $\br^d$ by $\bz^d$?

 As explained in the introduction, while this is known to be true for dimension $d=1$ or $d=2$, there are some 5 by 5 matrices for which such a $\mathcal D$ does not exist.

\begin{definition}
(a) We say that a subset $F\subset \br^d$ tiles $\br^d$ by a lattice $\Gamma$ iff the following two properties hold:
$$\br^d=\bigcup_{\gamma\in\Gamma}(F+\gamma)$$
$$(F+\gamma)\cap(F+\gamma')=\emptyset,\quad\gamma,\gamma'\in\Gamma,\gamma\neq\gamma'.$$
We say that $F$ tiles $\br^d$ up to measure zero by the lattice $\Gamma$ if these two properties hold up to Lebesgue measure zero.

(b) By a lattice we mean a rank-$d$ subgroup of $\br^d$. We shall be interested in sublattices $\Gamma\subset\bz^d$. For a fixed $\Gamma\subset\bz^d$ we say that $\Gamma$ is of index $k$ if the order of the quotient $\bz^d/\Gamma$ is $k$.
\end{definition}

\begin{lemma}\label{lem7_2}
Suppose a relatively compact subset $F\subset \br^d$ tiles by some lattice $\Gamma\subset\bz^d$. Then the lattice $\Gamma$ is of index $k$ if and only if the mapping $F\ni x\mapsto e^{2\pi ix}\in\bt^d$ is $k$-to-$1$ (up to measure zero).
\end{lemma}

\begin{proof}
It follows from the definition that $F$ tiles by $\Gamma$ iff the restriction to $F$ of the quotient mapping $\br^d\rightarrow\br^d/\Gamma$ is bijective up to measure zero. Hence the assertion that the given map is $k$-to-$1$ is equivalent to the natural mapping $\br^d/\Gamma\rightarrow\br^d/\bz^d$ being a $k$-fold cover; but this is so by the induced isomorphism $(\br^d/\Gamma)/(\br^d/\bz^d)\cong \bz^d/\Gamma$.
\end{proof}

\begin{proposition}\label{propmfd}
Let $\Omega:=\mathcal D^{\bn}$. Define the map $\mathfrak d:\br^d\times\Omega\rightarrow\mathcal S_A$ by
\begin{equation}\label{eqmathfrakd}
\mathfrak d(x,\omega_0\omega_1\dots)=(e^{2\pi ix},e^{2\pi i\tau_{\omega_0}x},e^{2\pi i\tau_{\omega_1}\tau_{\omega_0}x},\dots),\quad(x\in\br^d,\omega_0\omega_1\dots\in\Omega).
\end{equation}
\begin{enumerate}
	\item For each $x\in \br^d$, $k\in\bz^d$, and $\omega\in\Omega$ there is a unique $\omega'=\omega'(x,k,\omega)\in\Omega$ such that $\mathfrak d(x,\omega)=\mathfrak d(x+k,\omega')$. Moreover, 
	if $x,x'\in\br^d$ and $\omega,\omega'\in\Omega$, such that $\mathfrak d(x,\omega)=\mathfrak d(x',\omega')$, then 
	$x'=x+k$ for some $k\in\bz^d$ and $\omega'=\omega'(x,k,\omega)$.
	\item Let $F$ be a subset of $\br^d$. The restriction of the map $\mathfrak d$ to $F\times\Omega$ is injective if and only if $F\cap (F+k)=\emptyset$ for all $k\in\bz^d$, $k\neq0$.
	\item The restriction of $\mathfrak d$ to $F\times\Omega$ is onto if and only if 
	$$\bigcup_{k\in\bz^d}(F+k)=\br^d.$$
	\item The restriction of $\mathfrak d$ to $F\times\Omega$ is bijective if and only if $F$ tiles $\br^d$ by $\bz^d$-translations. 
	\item Define the map $\rho:\br^d\times\Omega\rightarrow\br^d\times\Omega$ 
	\begin{equation}\label{eqrho}
\rho(x,\omega_0\omega_1\dots)=(\tau_{\omega_0}x,\omega_1\omega_2\dots),\quad(x\in\br^d,\omega_0\omega_1\dots\in\Omega).
\end{equation}
Then 
\begin{equation}\label{eqrhosigma}
\mathfrak d\circ\rho=\sigma^{-1}\circ\mathfrak d.
\end{equation}
\end{enumerate}
\end{proposition}

\begin{proof}
(i) We want $\tau_{\omega_0}x\equiv\tau_{\omega_0'}(x+k)$. So $\at^{-1}\omega_0\equiv\at^{-1}(k+\omega_0')$, i.e., $\omega_0'\equiv \omega_0-k\mod A^T\bz^d$. Since $\mathcal D$ is a complete set of representatives for $\bz^d/A^T\bz^d$, there is a unique $\omega_0'\in\mathcal D$ such that this is satisfied. Proceeding by induction we see that $\omega_2',\omega_3',\dots$ can be uniquely constructed such that $e^{2\pi i\tau_{\omega_n}\dots\tau_{\omega_0}x}=e^{2\pi i\tau_{\omega_n'}\dots\tau_{\omega_0'}(x+k)}$ for all $n\in\bn$.

If $\mathfrak d(x,\omega)=\mathfrak d(x',\omega')$ then $e^{2\pi ix}=e^{2\pi ix'}$ so $x'=x+k$ for some $k\in\bz^d$. The rest follows from the uniqueness of $\omega'(x,k,\omega)$.

(ii) Suppose $\mathfrak d$ restricted to $F\times\Omega$ is injective. Then, by (i), we cannot have $x$ and $x+k$ in $F$ for some $k\neq 0$. Conversely, if $\mathfrak d$ is not injective on this set, then $\mathfrak d(x,\omega)=\mathfrak d(x',\omega')$ for some $x,x'\in F$ and $\omega,\omega'\in\Omega$. Using (i) again we get that $x'=x+l$ for some $l\in\bz^d$ so $F\cap (F+l)\neq\emptyset$. 

(iii) Suppose the restriction of $\mathfrak d$ to $F\times\Omega$ is onto. Then for all $y\in\br^d$, there is $x\in F$ and $\omega\in\Omega$ such that $\mathfrak d(x,\omega)=\hat i(y)$. Then $e^{2\pi ix}=e^{2\pi iy}$ so $y=x+k$ for some $k\in\bz^d$, and therefore $y\in F+k$. This shows that $\cup(F+k)=\br^d$. 

Conversely, take $(z_n)_{n\in\bn}\in\sa$. There exist $x_n\in\br^d$ such that $z_n=e^{2\pi ix_n}$ for all $n$. By hypothesis we can take $x_0\in F$. Since $z_1^A=z_0$ we have that, $A^Tx_1\equiv x_0$ so for some $\omega_0\in\mathcal D$, $\tau_{\omega_0}x_0\equiv x_1$. Then, by induction we can construct $\omega_n\in\mathcal D$ such that $\tau_{\omega_n}\dots\tau_{\omega_0}x_0\equiv x_{n+1}$. This proves that $\mathfrak d(x,\omega)=(e^{2\pi ix_n})_n=(z_n)_n$.

(iv) follows directly from (ii) and (iii). 

(v) requires nothing more that a simple computation.
\end{proof}

\begin{proposition}
Suppose $F$ is a subset of $\br^d$ that tiles $\br^d$ by a sublattice $\Gamma$ of $\bz^d$ with $|\bz^d/\Gamma|=N$. Then the restriction of the map $\mathfrak d$ to $F\times\Omega$ is $N$-to-$1$.
\end{proposition}

\begin{proof}
Since $\cup_{k\in\bz^d}(F+k)\supset \cup_{\gamma\in \Gamma}(F+\gamma)=\br^d$, it follows from Proposition \ref{propmfd}(iii) that the map is onto.

We claim that for each $x\in F$ there are exactly $N$ points $k\in\bz^d$ such that $x+k\in F$. Indeed, let $d_1,\dots d_{N}$ be a complete list of representatives for $\bz^d/\Gamma$. Then for each $i\in\{1,\dots,N\}$ there is a unique $\gamma_i\in\Gamma$ such that $x+d_i\in F+\gamma_i$. Then we can not have $d_i-\gamma_i=d_{i'}-\gamma_{i'}$ for $i\neq i'$ (that would imply $d_i\equiv d_{i'}\mod\Gamma$) so the points $d_i-\gamma_i$ are distinct and $x+(d_i-\gamma_i)\in F$. 

Now we can use Proposition \ref{propmfd}(i) to see that $\mathfrak d$ restricted to $F\times\Omega$ is $N$-to-$1$.
\end{proof}

Consider now the compact attractor $X(A^T,\mathcal D)$ of the iterated function system $(\tau_d)_{d\in\mathcal D}$ given in \eqref{eqxa}. It is known (see \cite{LaWa96c, LaWa97}) that $X(A^T,\mathcal D)$ always tiles $\br^d$ (up to measure zero) by some sublattice $\Gamma$ of $\bz^d$. 
\par The connection between lattice tilings and spectral theory was studied
systematically in \cite{Fug74} and \cite{Ped96}, and we will introduce spectrum in
subsection 6.2 below. From the choice of digit set $\mathcal D$ for a fixed matrix $A$,
we conclude that $X(A^T,\mathcal D)$ has non-empty interior. In fact the $d$-dimensional
Lebesgue measure of $X(A^T,\mathcal D)$ must be an integer. It is 1 if and only if $X(A^T,\mathcal D)$
tiles $\br^d$ by the "unit-lattice" $\bz^d$. By spectral theory we are referring to
the Hilbert space $L^2(X(A^T,\mathcal D))$.

\begin{definition}\label{deftile}
We say that $(A,\mathcal D)$ satisfy the tiling condition if $X(A^T,\mathcal D)$ tiles $\br^d$ (up to measure zero) by the lattice $\bz^d$.
\end{definition}

 In subsection 6.2 below we give examples for $d = 2$ of pairs $(A,\mathcal D)$ which do not satisfy the tiling condition.
      Nonetheless, even if some particular pair $(A,\mathcal D)$ in the plane does not satisfy the tiling condition, it will be possible to change the digit set $\mathcal D$ into a different one $\mathcal D'$, while keeping the matrix $A$ fixed, such that the modified pair $(A,\mathcal D')$ will satisfy the tiling condition. But by going to higher dimensions ($d = 5$) as we noted there are matrices $A$ for which no $\mathcal D$ may be chosen making $(A,\mathcal D)$ satisfy the tiling condition. 

\begin{lemma}\label{lemnooverlap}
For all $d,d'\in\mathcal D$, $d\neq d'$ the intersection $\tau_d(X(A^T,\mathcal D))\cap\tau_{d'}(X(A^T,\mathcal D))$ has Lebesgue measure zero.
\end{lemma}

\begin{proof}
We have the following relations, with $\mu$ the Lebesgue measure:
$$X(A^T,\mathcal D)=\bigcup_{d\in\mathcal D}\tau_d(X(A^T,\mathcal D)),$$
$$\mu(\tau_d(X(A^T,\mathcal D)))=\frac{1}{|\det A|}\mu(X(A^T,\mathcal D)).$$
As a result we get 
$$\mu(X(A^T,\mathcal D))=\sum_{d\in\mathcal D}\frac{1}{|\det A|}\mu(X(A^T,\mathcal D))-\mu(\mbox{ combined overlap sets}).$$
Therefore the combined overlap sets must have measure zero.
\end{proof}

\begin{proposition}\label{proprho}
Suppose $(A,\mathcal D)$ satisfy the tiling condition. Then the function $\rho$ defined in \eqref{eqrho} maps $X(A^T,\mathcal D)\times\Omega$ onto itself and the restriction of $\rho$ to $X(A^T,\mathcal D)\times\Omega$ is injective a.e. in the sense that the set of points $x\in X(A^T,\mathcal D)$ with the property that there exist $x'\in X(A^T,\mathcal D)$, $\omega,\omega'\in\Omega$  $\rho(x,\omega)=\rho(x',\omega')$, has Lebesgue measure zero. The inverse $\rho^{-1}$ of this restriction is defined by 
\begin{equation}\label{eqrhoi}
\rho^{-1}(x,\omega_0\omega_1\dots)=(A^Tx-\omega_x,\omega_x\omega_0\omega_1\dots),\quad(x\in X(A^T,\mathcal D),\omega_0\omega_1\dots\in\Omega),
\end{equation}
where $\omega_x$ is the unique element of $\mathcal D$ with the property $x\in\tau_{\omega_x}(X(A^T,\mathcal D))$.
\end{proposition}

\begin{proof}
Since \begin{equation}\label{eqatt}
X(A^T,\mathcal D)=\cup_{d\in\mathcal D}\tau_d(X(A^T,\mathcal D)),
\end{equation}
it follows that $\rho$ maps $X(A^T,\mathcal D)$ onto itself. 

Suppose now $\rho(x,\omega)=\rho(x',\omega')$ for $x,x'\in X(A^T,\mathcal D)$, $\omega,\omega'\in\Omega$, and $(x,\omega)\neq (x',\omega')$. So either $x\neq x'$ or $\omega\neq \omega'$.
When $\omega\neq \omega'$, since $\rho(x,\omega)=\rho(x',\omega')$ it follows that $\omega_1\omega_2\dots=\omega_1'\omega_2'\dots$ so $\omega_0\neq \omega_0'$. Also $\tau_{\omega_0}x=\tau_{\omega_0'}x'$. But $\tau_{\omega_0}(X(A^T,\mathcal D))\cap\tau_{\omega_0'}(X(A^T,\mathcal D))$ has measure zero (see Lemma \ref{lemnooverlap}), so $x$ must be in a set of measure zero. If $\omega=\omega'$ then $\omega_0=\omega_0'$ so $\tau_{\omega_0}x=\tau_{\omega_0'}x'$ implies $x=x'$. This proves the injectivity of $\rho$.

Since $\tau_d(X(A^T,\mathcal D))$ are mutually disjoint, the element $\omega_x$ of $\mathcal D$ is well defined. Then it is easy to check that $\rho^{-1}$ is indeed the inverse of the restriction of $\rho$.
\end{proof}

\subsection{Encodings of cyclic paths}
Let $\mathcal D$ be a complete set of representatives of $\bz^d/A^T\bz^d$. And let $\Omega:=\mathcal D^{\bn}$.

Consider now a cycle $C:=\{\zeta_0,\zeta_1,\dots,\zeta_{p-1}\}$ and suppose $\zeta_j=e^{2\pi i\theta_j}$ for some $\theta_j\in\br^d$. 
If in addition, $A$ and $\mathcal D$ satisfy the tiling condition we can pick $\theta_j$ in $X(A^T,\mathcal D)$. Then, since $z_1^A=z_0$, there is a $l_0\in\mathcal D$ such that $\theta_1=\tau_{l_0}\theta_0$. Continuing this process, we can find $l_0,l_1,\dots,l_{p-1}$ such that $\tau_{l_j}\theta_j=\theta_{j+1}$ for $j\in\{0,\dots,p-2\}$ and $\tau_{l_{p-1}}\theta_{p-1}=\theta_0$. Thus the point $\chi_C=(\zeta_0,\dots,\zeta_{p-1},\zeta_0,\dots,\zeta_{p-1},\zeta_0,\dots)$ from $\sa$ can be encoded by $(\theta_0,l_0\dots l_{p-1}l_0\dots l_{p-1},l_0\dots)$, i.e., by an infinite repetition of the finite word 
$l_0\dots l_{p-1}$.

\begin{definition}\label{defcycle}
A finite set in $\br^d$, $C=\{\theta_0,\dots,\theta_{p-1}\}$ is called a {\it cycle} if there exist $l_0,\dots,l_{p-1}\in\mathcal {D}$ such that $\tau_{l_0}\theta_0=\theta_{1},\tau_{l_1}\theta_{1}=\theta_{2},\dots,\tau_{l_{p-2}}\theta_{p-2}=\theta_{p-1}$ and $\tau_{l_{p-1}}\theta_{p-1}=\theta_{0}$. Thus $\theta_0$ is the fixed point of $\tau_{l_{p-1}}\dots\tau_{l_0}$, $\theta_1$ is the fixed point of $\tau_{l_0}\tau_{l_{p-1}}\dots\tau_{l_1}$, $\dots$, $\theta_{p-1}$ is the fixed point of $\tau_{l_{p-2}}\dots\tau_{l_0}\tau_{l_{p-1}}$. 

The points $\theta_0,\dots,\theta_{p-1}$ are called {\it cyclic points}. We say that $\theta_0$ is the {\it cyclic point associated to} $l_0\dots l_{p-1}$, and we say that $C=\{\theta_0,\dots,\theta_{p-1}\}$ is the {\it cycle} associated to $l_0\dots l_{p-1}$. 
\end{definition}

Let $C:=\{e^{2\pi i\theta_0},\dots,e^{2\pi i\theta_{p-1}}\}$ be a cycle associated to $l_0\dots l_{p-1}$. Take now a point $(z_n)_{n\in\bn}$ in $\sa(C)$. We want to see how the encodings of points in $\sa(C)$ look like.

By Theorem \ref{prop3_9}, the point $(z_n)_{n\in\bn}$ is in one of the sets $\sigma_A^{-j}(\chi_C)\hat i(\br^d)$. 
Suppose $z_0=e^{2\pi ix}$ for some $x\in\br^d$. Then $(z_n)_{n\in\bn}=\sigma_A^{-j}(\chi_C)\hat i(y)$ for some $y\in\br^d$, and
looking at the 0 position, $x\equiv \theta_j+y$ so $y=x+k-\theta_j$ for some $k\in\bz^d$. Thus
$$(z_n)_{n\in\bn}=\sigma_A^{-j}(\chi_C)\hat i(x+k-\theta_j)=\hat i_C(x+k-\theta_j,j).$$

On the other hand, according to the previous discussion, $(z_n)_{n\in\bn}$ is equal to $(e^{2\pi ix},e^{2\pi i\tau_{\omega_0}x},e^{2\pi i\tau_{\omega_1}\tau_{\omega_0}x},\dots)$ for some infinite word $\omega_0\omega_1\dots$. 
Thus we must have some precise correspondence between the pair $(k,j)\in\bz^d\times\bz_p$ and the infinite word $\omega_0\omega_1\dots\in\Omega$. Since $z_n$ is approaching the cycle $C$ as $n\rightarrow\infty$ one might expect that the infinite word $\omega_0\omega_1\dots$ ends in a repetition of the finite word $l_0\dots l_{p-1}$ that generates the cycle $C$. While this is often true, there might be some other cycles $C'$ that are congruent $\mod\bz^d$ to $C$, that will affect this encoding $\omega$. In any case, $\omega_0\omega_1\dots$ that corresponds to $(k,j)$ will be eventually periodic, and it will end in an infinite repetition of a finite word that corresponds to such a cycle $C'$.

\begin{definition}
We denote by $\underline{l_0\dots l_{p-1}}$ the infinite word in $\Omega$ obtained by the infinite repetition of the word $l_0\dots l_{p-1}$. Let
$$\Omega_C:=\{\omega_0\dots\omega_n\underline{l_0\dots l_{p-1}}\,|\,\omega_0,\dots,\omega_n\in\mathcal D,n\in\bn\},$$
i.e., the set of infinite words that end in an infinite repetition of the word $l_0\dots l_{p-1}$.
\end{definition}

There are some cycles which have cycle points that differ by integers. Such cycles would make our encoding ambiguous, so we avoid this situation. 
\begin{example} Let $d=1$, $A=2$ and $\mathcal D=\{0,3\}$. Then $\tau_0x=x/2$, $\tau_3=(x+3)/2$. Then it is easy to check that the attractor $X(A^T,\mathcal D)$ is $[0,3]$. The set $\{1,2\}$ is a cycle that corresponds to $\underline{30}$, and its points differ by an integer.
\end{example}
 
\begin{definition}
We say that the cycle $C=\{\theta_0,\dots,\theta_{p-1}\}$ is {\it simple} if $\theta_j\not\equiv\theta_{j'}\mod\bz^d$ for $j\neq j'$.
\end{definition}

Following \cite{BrJo99}, for a simple cycle $C=\{\theta_0,\dots,\theta_{p-1}\}$, we define an automorphism $\rc$ on the set $\bz^d-C$. Note that since the cycle is simple, the sets $\bz^d-\theta_j$ are mutually disjoint. The map $\rc$ is an extension of the division with remainder. Here we ``divide'' by $A^T$. For each point in $a-\theta_j\in\bz^d-\theta_j$, there is a unique ``quotient'' $b-\theta_{j+1}$ in $\bz^d-\theta_{j+1}$ and a unique ``remainder'' $d_0\in\mathcal D$ such that
$$a-\theta_j=A^T(b-\theta_{j+1})+d_0.$$
Then $\rc(a-\theta_j)$ is defined as the quotient $\rc(a-\theta_j)=b-\theta_{j+1}$.

We used here the fact that $A^T\theta_j\equiv\theta_{j-1}\mod\bz^d$, because $\tau_{l_{j-1}}\theta_{j-1}=\theta_j$, for all $j\in\bz$. Recall also that we use the notation $\theta_j:=\theta_{j\mod p}$ for $j\in\bz$.
\begin{definition}
Let $C=\{\theta_0,\dots,\theta_{p-1}\}$ be a simple cycle. On $\bz^d-C=\bigcup_{j=0}^{p-1}(\bz^d-\theta_j)$ we define the map $\mathcal R_C$ as follows: for each $a\in\bz^d$ and $j\in\{0,\dots,p-1\}$ there exist a unique $b\in\bz^d$ and $d_0\in\mathcal D$ such that
\begin{equation}\label{eqdefrc}
a-\theta_j=A^T(b-\theta_{j+1})+d_0.\mbox{ We define }\mathcal R_C(a-\theta_j):=b-\theta_{j+1}.
\end{equation}
Therefore $\mathcal R_C(a-\theta_j)$ is defined by
$$(a-\theta_j)-A^T\mathcal R_C(a-\theta_j)\in\mathcal D.$$
Also, $-\mathcal R_C(a-\theta_j)=\tau_{d_0}(-(a-\theta_j))$, where $d_0$ is the unique element of $\mathcal D$ such that
$\tau_{d_0}(-(a-\theta_j))\in\bz^d-\theta_{j+1}$. 
\end{definition}
 
The encoding of $(k,j)\in\bz^d\times\bz_p$ is obtained by a generalized Euclidean algorithm: take $k-\theta_j$ in $\bz^d-\theta_j$, then ``divide'' by $A^T$ and keep the remainder: $k-\theta_j=A^T\rc(k-\theta_j)+\omega_0$. Then take the quotient $\rc(k-\theta_j)$, divide by $A^T$ and keep the remainder $\omega_1$, and so on to infinity. The infinite sequence of remainders will give us $\omega$.

But first, we need some properties of the map $\rc$.
\begin{proposition}\label{proprc}
Let $C=\{\theta_0,\dots,\theta_{p-1}\}$ be a simple cycle. Let $X(A^T,\mathcal D)$ be the attractor of the iterated function system $(\tau_{d})_{d\in \mathcal D}$. 
\begin{enumerate}
\item A point $t\in C-\bz^d$ is a cycle point for the iterated function system $(\tau_{d})_{d\in\mathcal D}$ if and only if there is some $n\geq 1$ such that $\mathcal R_C^n(-t)=-t$, i.e., $-t$ is a periodic point for $\mathcal R_C$. Moreover if $t$ is associated to $m_0\dots m_{q-1}$ then $q$ is a multiple of $p$ and 
$$m_n=\rc^n(-t)-A^T\rc^{n+1}(-t),\quad(n\in\bn).$$
\item For every $t\in C-\bz^d$ there exists a $l\geq0$ such that $\mathcal R_C^l(-t)$ is periodic for $\mathcal R_C$, i.e., every point in $\bz^d-C$ is eventually periodic for $\mathcal R_C$. Moreover $-\mathcal R_C^l(-t)$ is in $(C-\bz^d)\cap X(A^T,\mathcal D)$.
\item 
The intersection $(C-\bz^d)\cap X(A^T,\mathcal D)$ consists exactly of negative the periodic points for $\mathcal R_C$.
\end{enumerate}

\end{proposition}

\begin{proof}
(i) If $t_0=t\in\theta_j-\bz^d$ is a cyclic point for $(\tau_{d})_{d\in \mathcal D}$, then $\tau_{m_0}t_0=t_1$, $\tau_{m_1}t_1=t_2,\dots,\tau_{m_{q-1}}t_{q-1}=t_0$ for some $m_0,\dots, m_{q-1}\in\mathcal D$ and some $t_1,\dots, t_{q-1}\in\br^d$. Then $t_{q-1}=A^Tt_0-m_{q-1}$ so $t_{q-1}\in\theta_{j-1}-\bz^d$, (because $A^T\theta_j\equiv \theta_{j-1}$). By induction $t_l\in\theta_{j+l-q}-\bz^d$ for $l\in\{q-1,q-2,\dots,0\}$. 

Since $t_0\in\theta_j-\bz^d$ and also $t_0\in\theta_{j-q}-\bz^d$, as the cycle is simple, it follows that $q$ must be a multiple of $p$. 

We have $\tau_{m_0}t_0=t_1$ so $-t_0=A^T(-t_1)+m_0$. Also $-t_1\in \bz^d-\theta_{j+1}$. This means that $\rc(-t_0)=-t_1$ and $m_0=-t_0-A^T\rc(-t_0)$. By induction $\rc(-t_{n})=-t_{n+1}$ and $m_n=-t_n-A^T\rc(-t_n)$. This proves one direction.

For the converse, if $\mathcal R_C^q(-t_0)=-t_0$, for some $t_0\in\theta_j-\bz^d$, then for each $n$ there is some $m_n\in\mathcal D$ such that:
$\rc^n(-t_0)=A^T\rc^{n+1}(-t_0)+m_n$. Thus the sequence $\{m_n\}$ has period $q$, and $\tau_{m_n}(-\rc^n(-t_0))=-\rc^{n+1}(-t_0)$, which proves that $t_0$ is in the cycle $\{-(-t_0),-\rc(-t_0),\dots,-\rc^{q-1}(-t_0)\}$.

(ii) Since $A$ is expansive there is a norm on $\br^d$ such that for some $0<c<1$, $\|\at^{-1}x\|\leq c\|x\|$ for all $x\in\br^d$. Then if $R>\frac{c\max_{d\in\mathcal D}\|d\|}{1-c}$, 
$$\tau_{d}(B(0,R))\subset B(0,R).$$
Indeed $\|\tau_{d}x\|\leq c(\|x\|+\|d\|)<c(R+\|d\|)<R$ for all $x\in B(0,R)$ and all $d\in\mathcal{D}$. 
\par
Take now $a-\theta_j\in\bz^d$. Take some $R>\max\{\|a-\theta_j\|,\frac{c\max_{d\in\mathcal D}\|d\|}{1-c}\}$. Then note that $\mathcal R_C(a-\theta_j)=-\tau_{d}(-(a-\theta_j))$ for some $d\in\mathcal D$. Therefore $\mathcal R_C$ maps $B(0,R)\cap\cup_j(\bz^d-\theta_j)$ into itself. So $\{\mathcal R_C^n(a-\theta_j)\,|\,n\in\bn\}$ is a finite set. 
Therefore there exists some $n\in\bn$, and $q\geq 0$ such that $\mathcal R_C^n(a-\theta_j)=\mathcal R_C^{n+q}(a-\theta_j)$. Thus $\mathcal R_C^n(a-\theta_j)$ is periodic.

From (i) we have that $-\rc^n(a-\theta_j)$ is cyclic for $(\tau_{d})_{d\in\mathcal D}$. So $-\rc^n(a-\theta_j)$ is in the attractor $X(A^T,\mathcal D)$. 

(iii) From (i) and (ii) it is clear that the periodic points for $\mathcal R_C$ lie in $(\bz^d-C)\cap(-X(A^T,\mathcal D))$. For the other inclusion take $t_1\in(\bz^d-C)\cap (-X(A^T,\mathcal D))$. Then using the formula \eqref{eqxa}, there exist $d_1,d_2\dots\in\mathcal D$ such that
$$-t_1=\at^{-1}d_1+\at^{-2}d_2+\dots$$
Let $-t_n:=\at^{-1}d_{n}+\at^{-2}d_{n+1}+\dots$. We have 
\begin{equation}\label{eqtn}
A^T(t_n)+d_n=t_{n+1},\quad(n\in\bn).
\end{equation}
Since $t_1\in\bz^d-C$, equation \eqref{eqtn} implies that $t_2$ is in $\bz^d-C$ and $t_1=\mathcal R_C(t_2)$. By induction $t_{n+1}$ is in $\bz^d-C$ and $\rc(t_{n+1})=t_n$ for all $n\in\bn$. But we have also $t_n\in -X(A^T,\mathcal D)$. And, since $(\bz^d-C)\cap (-X(A^T,\mathcal D))$ is finite, there exist $n,m\geq 1$ such that $t_n=t_{n+m}$. This implies that 
$\mathcal R_C^m(t_n)=\mathcal R_C^m(t_{n+m})=t_n$. Since $\mathcal R_C^{n-1}(t_n)=t_1$, it follows that $t_1$ is periodic for $\mathcal R_C$.
\end{proof}

\begin{theorem}\label{thenccycl}
Let $C:=\{\theta_0,\dots,\theta_{p-1}\}$ be a simple cycle. 
\begin{enumerate}
\item 
For each $k\in\bz^d$ and each $j\in\bz_p$ there is a unique $\omega(k,j)=\omega_0\omega_1\dots\in\Omega$ such that for all $x\in\br^d$,
\begin{equation}\label{eqomegak}
(e^{2\pi i x},e^{2\pi i\tau_{\omega_0}x},e^{2\pi i\tau_{\omega_1}\tau_{\omega_0}x},\dots)=\hat i_C(x+k-\theta_j,j)=\sigma_A^{-j}(\chi_C)\hat i(x+k-\theta_j)=(e^{2\pi i(\at^{-n}(x+k-\theta_j)+\theta_{j+n})})_{n\in\bn}.
\end{equation}
Moreover there exists a cycle $C'\in (C-\bz^d)\cap X(A^T,\mathcal D)$ such that $\omega(k,j)\in\Omega_{C'}$. 
\item The infinite word $\omega(k,j)$ can be constructed as follows: 
\begin{equation}
	\omega_n=\mathcal R_C^n(k-\theta_j)-A^T\mathcal R_C^{n+1}(k-\theta_j),\quad(n\in\bn).
	\label{eq:defomegank}
\end{equation}
\item
Suppose $C'$ is a cycle in $(C-\bz^d)\cap X(A^T,\mathcal D)$, and $\omega\in\Omega_{C'}$. Then there is a unique $(k(\omega),j(\omega))\in\bz^d\times\bz_p$ such that for all $x\in\br^d$,
\begin{equation}\label{eqkomega}
(e^{2\pi i x},e^{2\pi i\tau_{\omega_0}x},e^{2\pi i\tau_{\omega_1}\tau_{\omega_0}x},\dots)=\hat i_C(x+k(\omega)-\theta_{j(\omega)},j(\omega)).
\end{equation}
\item 
$(k(\omega),j(\omega))$ can be constructed as follows: if $\omega\in\Omega_{C'}$ has the form $\omega_0\dots\omega_{np-1}\underline{m_0\dots m_{q-1}}$, then the fixed point $\eta_0$ of $\tau_{m_{q-1}}\dots\tau_{m_0}$ belongs to $\theta_{j(\omega)}-\bz^d$ for some unique $j(\omega)\in\bz_p$. 
And 
\begin{equation}\label{eqdefkomega}
k(\omega)=\omega_0+\dots+\at^{np-1}\omega_{np-1}+\theta_{j(\omega)}-\at^{np}\eta_0.
\end{equation}
\item Let
$$\tilde \Omega_C:=\bigcup\{\Omega_{C'}\,|\, C'\mbox{ cycle in }(C-\bz^p)\cap X(A^T,\mathcal D)\}.$$
The maps $$\mathfrak e_C:\bz^d\times\bz_p\rightarrow \tilde \Omega_C, \mathfrak e_C(k,j)=\omega(k,j),$$ and $$\mathfrak d_C:\tilde \Omega_C\rightarrow\bz^d\times\bz_p, \mathfrak d_C(\omega)=(k(\omega),j(\omega))$$ are inverse to each other.
\end{enumerate}
\end{theorem}

\begin{proof}
Let $(k,j)\in\bz^d\times\bz_p$ and let $\omega(k,j)$ be defined as in (ii). We prove that the relation \eqref{eqomegak} is satisfied. We have
$$A^T\mathcal R_C(k-\theta_j)+\omega_0=k-\theta_j, A^T\rc^2(k-\theta_j)+\omega_1=\rc(k-\theta_j),\dots$$
Therefore
$$\rc(k-\theta_j)=\at^{-1}(k-\theta_j)-\at^{-1}\omega_0, \rc^2(k-\theta_j)=\at^{-2}(k-\theta_j)-\at^{-2}\omega_0-\at^{-1}\omega_1,\dots$$
By induction 
$$\rc^n(k-\theta_j)=\at^{-n}(k-\theta_j)-\at^{-n}\omega_0-\dots-\at^{-1}\omega_{n-1}.$$
So
$$\tau_{\omega_{n-1}}\dots\tau_{\omega_0}x=\at^{-n}x
+\at^{-n}\omega_0+\dots+\at^{-1}\omega_{n-1}=\at^{-n}x+\at^{-n}(k-\theta_j)-\rc^n(k-\theta_j).$$
But $\rc^n(k-\theta_j)\in\bz^d-\theta_{j+n}$ so 
$$\tau_{\omega_{n-1}}\dots\tau_{\omega_0}x\equiv \at^{-n}(x+k-\theta_j)+\theta_{j+n}.$$
Therefore the relation \eqref{eqomegak} is satisfied. 

Next we prove the uniqueness of $\omega$. Suppose $\omega'\in\Omega$ also satisfies \eqref{eqomegak}. Then 
$\tau_{\omega_0'}x\equiv\tau_{\omega_0}x$ so $\at^{-1}\omega_0'\equiv\at^{-1}\omega_0$ which implies that 
$\omega_0'-\omega_0\in A^T\bz^d$. Since $\mathcal D$ is a complete set of representatives for $\bz^d/A^T\bz^d$, it follows that $\omega_0'=\omega_0$. By induction $\omega_n'=\omega_n$ so $\omega'=\omega$.

To see that $\omega(k,j)$ is in some $\omega_{C'}$ for a cycle $C'$ in $(C-\bz^d)\cap X(A^T,\mathcal D)$ we use Proposition \ref{proprc}. There exists an $l$ such that $\rc^l(k-\theta_j)$ is periodic for $\mathcal R_C$, so $-\mathcal R_C^l(k-\theta_j)$ is a cycle point for $(\tau_{d})_{d\in\mathcal D}$. Let $C'$ its corresponding cycle. By Proposition \ref{proprc}, $C'$ is contained in $(C-\bz^d)\cap X(A^T,\mathcal D)$. By the construction of $\omega(k,j)$ given in (ii), and by Proposition \ref{proprc}(i), we see that $\omega(k,p)\in\Omega_{C'}$.

Now let $\omega$ be of the form given in (iii). And let $(k,j):=(k(\omega),j(\omega))$ be as in (iv). Since $\eta_0\in \theta_{j}-\bz^d$ it follows that $\at^{np}\eta_0$ is also in $\theta_{j}-\bz^d$ so $k=k(\omega)$ is indeed an integer. 

We check that if for $(k,j)$ we construct $\nu=\nu_0\nu_1...=\omega(k,j)$ as in (ii), then $\nu=\omega$. This will prove also (v). From \eqref{eqdefkomega} we have with $t_0=-(k-\theta_j)$,
$$t_1:=\tau_{\omega_0}(t_0)=\tau_{\omega_0}(-(k-\theta_j))=-\omega_1-\dots-\at^{np-2}\omega_{np-1}+\at^{np-1}\eta_0\in\theta_{j+1}-\bz^d,$$
since $\eta_0\in \theta_j-\bz^d$ implies that $\at^{np-1}\eta_0\in\theta_{j-(np-1)}-\bz^d=\theta_{j+1}-\bz^d$. This shows that $\rc(-t_0)=-t_1$, and $\nu_0=\omega_0$. By induction we obtain $\nu_1=\omega_1,\dots,\nu_{np-1}=\omega_{np-1}$ and that $-\eta_0=\rc^{np}(-t_0)$. And since the cycle $C'$ of $\eta_0$ is given by $\underline{m_0\dots m_{q-1}}$, it follows by Proposition \ref{proprc}(i) that $\nu=\omega_0\dots\omega_{np-1}\underline{m_0\dots m_{q-1}}$.

For the uniqueness of $(k,j)$, suppose $(k,j)$ and $(k',j')$ satsify \eqref{eqkomega}. Then 
$$\at^{-n}(k-\theta_j)+\theta_{n+j}\equiv \at^{-n}(k'-\theta_{j'})+\theta_{n+j'},\quad(n\in\bn).$$
Taking the limit as $np\rightarrow\infty$ we obtain $\theta_j-\theta_{j'}\in\bz^d$. Since the cycle is simple, $j=j'$. 
Therefore
$$\at^{-n}(k-\theta_j)\equiv \at^{-n}(k'-\theta_{j}),\quad(n\in\bn).$$
But this means that $\hat i(k-\theta_j)=\hat i(k'-\theta_{j})$, and by Proposition \ref{propembed1} $\hat i$ is injective so $k=k'$. 
\end{proof}
 We summarize our results in the following corollary.
\begin{corollary}\label{corsum}
Suppose $(A,\mathcal D)$ satisfy the tiling condition (Definition \ref{deftile}). Let $C=\{\theta_0,\dots,\theta_{p-1}\}$ be a simple cycle. 
Let 
$$\tilde \Omega_C:=\bigcup\{\Omega_{C'}\,|\, C'\mbox{ cycle in }(C-\bz^d)\cap X(A^T,\mathcal D)\}.$$
\begin{enumerate}
\item
The maps
$$\mathfrak d:X(A^T,\mathcal D)\times\tilde \Omega_C\rightarrow\sa(C),\quad \mathfrak d(x,\omega)=(e^{2\pi ix},e^{2\pi i\tau_{\omega_0}x},e^{2\pi i\tau_{\omega_1}\tau_{\omega_0}x},\dots),$$
$$\hat i_C:\br^d\times\bz_p\rightarrow \sa(C),\quad \hat i_C(x,j)=(e^{2\pi i(\at^{-n}x+\theta_{n+j})})_{n\in\bn}$$
are bijections. 
\item 
$$\hat i_C^{-1}(\mathfrak d(x,\omega))=(x-\theta_{j(\omega)}+k(\omega),j(\omega)),\quad(x\in X(A^T,\mathcal D),\omega\in\tilde\Omega_C),$$
where $k(\omega),j(\omega)$ are defined in Theorem \ref{thenccycl}(iv). 
$$\mathfrak d^{-1}(\hat i_C(x,j))=(y,\omega(k,j)),\quad(x\in\br^d, j\in\bz_p),$$
where $y\in X(A^T,\mathcal D)$, $k\in\bz^d$ are uniquely defined by $x+\theta_j=y+k$, and $\omega(k,j)$ is defined in Theorem \ref{thenccycl}(ii).
\item The following diagram is commutative:
$$
\begin{CD}
X(A^T,\mathcal D)\times\tilde \Omega_C @>\mathfrak d>>  \sa(C) @<\hat i_C<< \br^d\times\bz_p\\
@VV\rho^{-1}V 		@VV\sigma_A V @VV \alpha_{A,p} V\\
X(A^T,\mathcal D)\times\tilde \Omega_C@>\mathfrak d>> \sa(C) @<\hat i_C<< \br^d\times\bz_p
\end{CD}
$$
\end{enumerate}
\end{corollary}
\begin{corollary}
With the notations in Theorem \ref{thenccycl}, we have
\begin{equation}\label{eqjomega}
j(\omega_1\omega_2\dots)=j(\omega_0\omega_1\dots)+1,\quad(\omega_0\omega_1\dots\in\tilde \Omega_C).
\end{equation}
\begin{equation}
\at^{-n}(x-\theta_{j(\omega_0\omega_1\dots)}+k(\omega_0\omega_1\dots))=\tau_{\omega_{n-1}}\dots\tau_{\omega_0}x-\theta_{j(\omega_0\omega_1\dots)+n}+k(\omega_n\omega_{n+1}\dots)
\end{equation}
for all $x\in X(A^T,\mathcal D),\omega_0\omega_1\dots\in\tilde\Omega_C$.
\end{corollary}

\begin{proof}
We apply the commutative diagram in Corollary \ref{corsum} to $\rho^n$:
$$\hat i_C^{-1}\mathfrak d\rho^n(x,\omega)=\hat i_C\mathfrak d(\tau_{\omega_{n-1}}\dots\tau_{\omega_0}x,\omega_n\omega_{n+1}\dots)=(\tau_{\omega_{n-1}}\dots\tau_{\omega_0}x-\theta_{j(\omega_n\omega_{n+1}\dots)}+k(\omega_n\omega_{n+1}\dots),j(\omega_n\omega_{n+1}\dots)).$$
$$\alpha_{A,p}^{-n}\hat i_C^{-1}\mathfrak d(x,\omega)=\alpha_{A,p}^{-n}(x-\theta_{j(\omega_0\omega_1\dots)}+k(\omega_0\omega_1\dots),j(\omega_0\omega_1\dots))=$$$$
(\at^{-n}(x-\theta_{j(\omega_0\omega_1\dots)}+k(\omega_0\omega_1\dots)),j(\omega_0\omega_1\dots)+n).$$
Since the two quantities are equal to each other (by the commutative diagram in Corollary \ref{corsum}), the relations follow.
\end{proof}

 With the aid of our cycles and associated encoding/decoding
mappings we are now able to state our main result regarding super
representations. Notice that the introduction of cycles yields the
following improvement of Theorem \ref{prop3_9} in section 4 above.

\begin{corollary}
Suppose $(A,\mathcal D)$ satisfies the tiling condition, and let $C$ be a simple cycle of length $p$. On $X(A^T,\mathcal D)\times\tilde\Omega_C$ define the measure $\breve\mu$ by
$$\int_{X(A^T,\mathcal D)\times\tilde\Omega_C}f\,d\breve\mu=\int_{X(A^T,\mathcal D)}\sum_{\omega\in\tilde\Omega_C}f(x,\omega)\,dx.$$
Define the operators $\breve T_k$, $k\in\bz^d$ and $\breve U$ on $L^2(X(A^T,\mathcal D)\times\tilde\Omega_C,\breve\mu)$ by
$$\breve T_kf(x,\omega)=e^{2\pi ik\cdot x}f(x,\omega),\quad(x\in X(A^T,\mathcal D),\omega\in\tilde\Omega_C),$$
$$\breve Uf=\sqrt{|\det A|}f\circ\rho^{-1}.$$
Then $\{\breve T_k,\breve U\}$ define a unitary representation of $G_A$ and $\mathcal W:L^2(\br^d\times\bz_p)\rightarrow L^2(X(A^T,\mathcal D),\breve\mu)$, $\mathcal Wf=f\circ\mathfrak d\circ\hat i_C^{-1}$ is an isomorphism that intertwines this representation with the one in Definition \ref{defhc}.
\end{corollary}

 When $(A,\mathcal D)$ satisfy the tiling condition we can say a bit more about the possible extra cycles in $(C-\bz^d)\cap X(A^T,\mathcal D)$:
 
 \begin{proposition}
 Suppose $(A,\mathcal D)$ satisfies the tiling condition. Assume that there is a cycle point $\theta_0\in X(A^T,\mathcal D)$ such that $\theta_0-k\in X(A^T,\mathcal D)$ for some $k\in\bz^d$, $k\neq 0$. Then the entire cycle of $\theta_0$ is on the boundary of $X(A^T,\mathcal D)$.
 \end{proposition}
 
 \begin{proof}
 Let $X(A^T,\mathcal D)^\circ$ denote the interior of $X(A^T,\mathcal D)$. 
 We will prove first that if a point $x$ is in $X(A^T,\mathcal D)\cap(X(A^T,\mathcal D)+k)$ with $k\in\bz^d$, $k\neq0$, then $x$ is on the boundary of $X(A^T,\mathcal D)$. Suppose not, then $x\in X(A^T,\mathcal D)^\circ$. By \cite{LaWa96c} we know that the closure of $X(A^T,\mathcal D)^\circ$ is $X(A^T,\mathcal D)$. This implies that the neighborhood $X(A^T,\mathcal D)^\circ$ of $x$ must intersect the set $X(A^T,\mathcal D)^\circ +k$. But since $X(A^T,\mathcal D)$ tiles $\br^d$ by $\bz^d$, this implies that the interiors of $X(A^T,\mathcal D)$ and $X(A^T,\mathcal D)+k$ cannot intersect (the intersection would have positive Lebesgue measure).  So $x$ must be on the boundary of $X(A^T,\mathcal D)$. 
 
 Now consider $\theta_0$ and let $C=\{\theta_0,\dots,\theta_{p-1}\}$ be the cycle of $\theta_0$ and let $l_0,\dots l_{p-1},$ the corresponding  digits. We have $\tau_{l_{p-1}}\dots\tau_{l_j}\theta_j=\theta_0$ for all $j\in\{0,\dots,p-1\}$. 
 Suppose one of the points $\theta_j$ of the cycle $C$ is in $X(A^T,\mathcal D)^\circ$. Since
 $$X(A^T,\mathcal D)=\bigcup_{d\in\mathcal D}\tau_d(X(A^T,\mathcal D)),$$
 we obtain that $\tau_d(X(A^T,\mathcal D)^\circ)\subset X(A^T,\mathcal D)^\circ$ ($\tau_d$ is a homeomorphism). 
 So if $\theta_j$ is an interior point for $X(A^T,\mathcal D)$, then $\theta_0=\tau_{l_{p-1}}\dots \tau_{l_j}\theta_j$ is also in the interior of $X(A^T,\mathcal D)$. This contradiction implies that $C$ is contained in the boundary of $X(A^T,\mathcal D)$. 
 \end{proof}

\par To help the reader appreciate our encoding results we present some examples which at the same time stress tiles versus spectrum. Since the technical points are illustrated already for the real line we begin with dimension one, and then turn to the plane $\br^2$.

\begin{example}
Let us take $d=1$, $A=2$ and $\mathcal D=\{0,1\}$. The maps are $\tau_0x=x/2$, $\tau_1x=(x+1)/2$. Consider the simple cycle $C:=\{\theta_0\}=\{0\}$. It corresponds to $\underline 0$. 

The attractor $X(A^T,\mathcal D)$ is $[0,1]$. The intersection $(C-\bz)\cap X(A^T,\mathcal D)=(-\bz)\cap[0,1]=\{0,1\}$, so it consists of the cycles $C'$: $\{0\}$ and $\{1\}$, which correspond to $\underline 0$ and $\underline 1$ respectively. Therefore 
$\tilde\Omega_0=\Omega_0\cup\Omega_1$, i.e., the words that end in an infinite repetition of $0$ or an infinite repetition of $1$. 

We have the map $$\hat i_0^{-1}\mathfrak d: [0,1)\times(\Omega_0\cup\Omega_1)\rightarrow \br,\quad \hat i_0^{-1}\mathfrak d(x,\omega)=x+k(\omega),$$
(we used $[0,1)$ here instead of $[0,1]$ to have that the map $\hat i_0^{-1}\mathfrak d$ is a true bijection, not just up to measure $0$),
and with formula \eqref{eqdefkomega}:
$$k(\omega_0\dots\omega_n\underline0)=\omega_0+2\cdot\omega_1+\dots+2^n\omega_n,\quad k(\omega_0\dots\omega_n\underline1)=\omega_0+2\cdot\omega_1+\dots+2^n\omega_n-2^{n+1}.$$

Applying the commutative diagram in Corollary \ref{corsum} to $\rho^n$, we have $\hat i_0^{-1}\mathfrak d \rho^n(x,\omega)=\frac{1}{2^n}\hat i_0^{-1}\mathfrak d(x,\omega),$ which implies that 
$$\frac{1}{2^n}(x+k(\omega_0\omega_1\dots))=\tau_{\omega_{n-1}}\dots\tau_{\omega_0}x+k(\omega_n\omega_{n+1}\dots),\quad (x\in[0,1), \omega\in \Omega_0\cup \Omega_1).$$
 
\end{example}

\begin{example}
Let $d=1$, $A=2$ and $\mathcal D=\{0,1\}$. Consider the simple cycle $C:=\{\theta_0,\theta_1\}=\{1/3,2/3\}$. It corresponds to $\underline{10}$ (because $\tau_1(1/3)=2/3$, $\tau_0(2/3)=1/3$. The attractor is $X(A^T,\mathcal D)=[0,1]$. Then $(C-\bz)\cap [0,1]=C$, so 
$\tilde \Omega_C=\Omega_C$, i.e., the words that end in $\underline{10}$.

We have that the map $$\hat i_C^{-1}\mathfrak d:[0,1)\times\Omega_C\rightarrow \br\times \bz_2,\quad
\hat i_C^{-1}\mathfrak d(x,\omega)=(x-\theta_{j(\omega)}+k(\omega),j(\omega))$$
is a bijection, and 
$$j(\omega_0\dots\omega_{2n-1}\underline{10})=0,\quad j(\omega_0\dots\omega_{2n-1}\underline{01})=1,$$
$$k(\omega_0\dots\omega_{2n-1}\underline{10})=\omega_0+2\cdot\omega_1+\dots+2^{2n-1}\omega_{2n-1}+\frac13-2^{2n}\cdot\frac13,$$
$$k(\omega_0\dots\omega_{2n-1}\underline{01})=\omega_0+2\cdot\omega_1+\dots+2^{2n-1}\omega_{2n-1}+\frac23-2^{2n}\cdot\frac23.$$

The map $\Omega_C\ni \omega\mapsto (k(\omega),j(\omega))\in\bz\times\bz_2$ is a bijection. 

As an example, let us show how to compute the $\omega\in\Omega_C$ associated to $(k,j)=(15,0)$.
Take $k-\theta_j=15-\frac13$. We want a $k_1\in\bz$ and $\omega_0\in \mathcal D$ such that $15-\frac13=2(k_1-\frac23)+\omega_0$. We have 
\begin{align*}
15-\frac13&=2(8-\frac23)+0& 1-\frac13&=2(1-\frac23)+0\\
8-\frac23&=2(4-\frac13)+0& 1-\frac23&=2(0-\frac13)+1\\
4-\frac13&=2(2-\frac23)+1& 0-\frac13&=2(0-\frac23)+1\\
2-\frac23&=2(1-\frac13)+0& 0-\frac23&=2(0-\frac13)+0\\
& &\vdots
\end{align*}
\end{example}

Thus $\omega(15,0)=001001\underline{10}$.

\begin{remark}
Our next example is in the plane, but it illustrates a more general picture in $\br^d$ for any $d$. Start with a given pair $(A, \mathcal D)$ with the matrix $A$ assumed expansive, and $\mathcal D$ a chosen complete digit set, i.e., in bijective correspondence with the points in $\bz^d/A^T\bz^d$. So in particular,  $|\mathcal D| = |\det A|$. In general it is not true that the same set $\mathcal D$ is a digit set for $A$, i.e., that it is a bijective image of $\bz^d/A\bz^d$. 
 Here for $d = 2$, we give an explicit geometric representation of a pair $(A, \mathcal D)$ for which the same $\mathcal D$ is a digit set for both the radix representation with $A$ and with the transposed matrix $A^T$. Hence we get two attractors  $X(A^T,\mathcal D)$  and  $X(A,\mathcal D)$. Both will be referred to as Cloud Nine, a left-handed version, and a right-handed version. That is because there are nine integer points, i.e., the intersections with $\bz^2$ consists of nine points, and it is the same set for the two fractals. For each, there are three one-cycles, and one six-cycle. While the six-cycles (for $A$ and for $A^T$) are the same as sets, we will see that they are traveled differently under the actions discussed in our encodings from sections 5 and 6 from above; the difference being essentially a reversal of orientation. Hence our encoding with infinite words in letters from $\mathcal D$ will also be different for the two cases, and the details are worked out below. Recall the attractor  $X(A^T,\mathcal D)$ is the set of ``fractions'' for our digital representation of points in $\br^2$. The attractor $X(A^T,\mathcal D)$ is also an affine Iterated Function System (IFS) based on $(A,\mathcal D)$. Thus the Cloud Nine examples further illustrate the intricate part  played by the cycles in $\bz^2$ for the initial $(A,\mathcal D)$-IFS. In each case, using these cycles we are able to write down formulas for the two maps which do the encoding as well as the decoding in our positional $\mathcal D$-representation. 
\end{remark}

\par For use of matrices in radix representation, the distinction
between the radix matrix $A$ and its transpose $A^T$ is important. First the two
matrices sit on separate sides in a Fourier duality; and secondly, even if 
the chosen set of digits $\mathcal D$ is the same, the two attractors may be different.
In fact, in general the same $\mathcal D$ may not work for both $A$ and $A^T$. There is not
a natural connection between the two quotients $\bz^d/A\bz^d$  and $\bz^d/A^T\bz^d$
, i.e., the one for $A$ and the other one for the transposed.
\par But in the particular 2D example below, Example \ref{a0}, called
Cloud-Nine, one may check by hand that, for this matrix $A$, with $\det A=5$, 
each of the two quotients $\bz^2/A\bz^2$  and $\bz^2/A^T\bz^2$  are in bijective
correspondence with the same subset $\mathcal D$ in $\bz^2$. (See details!)
\par As a result it makes sense to analyze the two different Hutchinson
attractors $X(A^T, \mathcal D)$ and $X(A, \mathcal D)$; both compact and with non-empty interior.
The first one, in a different context was studied earlier in \cite{BrJo99} and
\cite{Jor03}, but both are interesting. Note that while two references \cite{BrJo99}
and \cite{Jor03} use these examples the questions addressed in these papers are
completely different.
\par As we see, there is an intriguing connection between cycles,
solenoids, and encodings for the two.

\begin{example}\label{a0}
Take $d=2$, $A=\left(\begin{array}{cc}1&-2\\2&1\end{array}\right)$ so $A^T=\left(\begin{array}{cc}1&2\\-2&1\end{array}\right)$, and $\mathcal D=\left\{\vectr00,\vectr{\pm3}{0},\vectr{0}{\pm2}\right\}$. 

We consider the trivial cycle $C:=\{\vectr00\}$. We want to compute $(C-\bz^d)\cap X(A^T,\mathcal D)$, i.e., $\bz^d\cap X(A^T,\mathcal D)$. 

First, we need to locate the attractor $X(A^T,\mathcal D)$. For this we use the proof of Proposition \ref{proprc}(ii), and conclude that if $R:=\frac{\|\at^{-1}\|\max_{d\in\mathcal D}\|d\|}{1-\|\at^{-1}\|}$ then the ball $B(0,R)$ is invariant under all the maps $\tau_d$ which implies that $X(A^T,\mathcal D)$ is contained in this ball.

 Since $\|\at^{-1}\|=\frac{1}{\sqrt{5}}$, we conclude that $R=\frac3{\sqrt{5}-1}=2.427\dots$. There are 21 points in $\bz^d\cap B(0,R)$ and we can check how $\rc$ acts on each of them. 

 If we want $\vectr{x}{y}=A^T\vectr{a}{b}+\vectr{d_1}{d_2}$ with $\vectr{a}{b}\in\bz^2$ and $\vectr{d_1}{d_2}\in\mathcal D$, then we must have
 \begin{equation}\label{eqcl9}
 \frac{(x-d_1)-2(y-d_2)}5=a,\quad \frac{2(x-d_1)+(y-d_2)}5=b.
 \end{equation}
 Thus, given $\vectr{x}{y}$, to find $\rc\vectr xy=\vectr ab$ and $\vectr{d_1}{d_2}\in\mathcal D$, we first look for an element in $\mathcal D$ with $d_1-2d_2\equiv x-2y\mod 5$ and $2d_1+d_2\equiv 2x+y\mod 5$. Then we compute $a,b$ as in \eqref{eqcl9}.
 
 It is interesting to note here also that if $d_1-2d_2\equiv x-2y\mod 5$ then the other equivalence $\mod 5$ is satisfied too. This is because $2d_1+d_2\equiv 2(d_1-2d_2)\mod 5$ etc.
 We can use the following table
 \begin{center}
\begin{tabular}{|c|c|c|c|c|c|}
\hline
$\vectr{d_1}{d_2}$&$\vectr{0}{0}$&$\vectr{3}{0}$&$\vectr{-3}{0}$&$\vectr{0}{2}$&$\vectr{0}{-2}$\\
\hline
$(d_1-2d_2)\mod 5$&0&3&2&1&4\\
\hline
\end{tabular}
\end{center}
 
 We apply these ideas to the points in $\bz^2\cap B(0,R)$:

 $\vectr00=A^T\vectr00+\vectr00$ so $\rc\vectr00=\vectr00$, and $-\vectr00$ is a cycle that corresponds to $\underline{\vectr00}$. 
 
 $\vectr10=A^T\vectr10+\vectr02$ so $\rc\vectr10=\vectr10$, and $\vectr{-1}0$ is a cycle that corresponds to $\underline{\vectr02}$.
 
 $\vectr{-1}0=A^T\vectr{-1}0+\vectr0{-2}$ so $\rc\vectr{-1}{0}=\vectr{-1}{0}$, and $\vectr{1}{0}$ is a cycle that corresponds to $\underline{\vectr{0}{-2}}$.
 
 $\vectr01=A^T\vectr{-1}{-1}+\vectr30$, $\vectr{-1}{-1}=A^T\vectr{1}{-1}+\vectr02$, $\vectr{1}{-1}=A^T\vectr{0}{-1}+\vectr30$, $\vectr{0}{-1}=A^T\vectr{1}{1}+\vectr{-3}{0}$, $\vectr{1}{1}=A^T\vectr{-1}{1}+\vectr{0}{-2}$, $\vectr{-1}{1}=A^T\vectr{0}{1}+\vectr{-3}{0}$. 
 So 
 $$\vectr01\stackrel{\mathcal R_C}{\rightarrow}\vectr{-1}{-1}\stackrel{\mathcal R_C}{\rightarrow}\vectr{1}{-1}\stackrel{\mathcal R_C}{\rightarrow}\vectr{0}{-1}\stackrel{\mathcal R_C}{\rightarrow}\vectr11\stackrel{\mathcal R_C}{\rightarrow}\vectr{-1}{1}\stackrel{\mathcal R_C}{\rightarrow}\vectr{0}{1}$$
 
 and $\left\{\vectr{0}{-1},\vectr{1}{1},\vectr{-1}{1},\vectr{0}{1},\vectr{-1}{-1},\vectr{1}{-1}\right\}$ is a cycle that corresponds to the word $$\underline{\vectr30\vectr02\vectr30\vectr{-3}0\vectr0{-2}\vectr{-3}{0}}.$$
 
 Similar computations show that $\vectr20\stackrel{\mathcal R_C}{\rightarrow}\vectr12\stackrel{\mathcal R_C}{\rightarrow}\vectr02\stackrel{\mathcal R_C}{\rightarrow}\vectr00$, $\vectr{-2}{0}\stackrel{\mathcal R_C}{\rightarrow}\vectr{-1}{-2}\stackrel{\mathcal R_C}{\rightarrow}\vectr{0}{-2}\stackrel{\mathcal R_C}{\rightarrow}\vectr00$, $\vectr{1}{-2}\stackrel{\mathcal R_C}{\rightarrow}\vectr10$, $\vectr21\stackrel{\mathcal R_C}{\rightarrow}\vectr01$, $\vectr{2}{-1}\stackrel{\mathcal R_C}{\rightarrow}\vectr{0}{1}$, $\vectr{-2}{1}\stackrel{\mathcal R_C}{\rightarrow}\vectr{0}{-1}$, $\vectr{-2}{-1}\stackrel{\mathcal R_C}{\rightarrow}\vectr{0}{-1}$.
 
 Thus we have three cycles of length one, and one cycle of length 6.

Consider now the pair $(A^T, \mathcal D)$, that is replace the matrix $A$ by $A^T$, and $\mathcal D$ is also a complete set of representatives for $\bz^d/\at^T\bz^d$. As above we get the same cycles:\\
$\vectr00$ is a one-cycle that corresponds to $\underline{\vectr00}$\\
$\vectr10$ is a one-cycle that corresponds to $\underline{\vectr0{2}}$\\
$\vectr{-1}0$ is a one-cycle that corresponds to $\underline{\vectr{0}{-2}}$\\
We also obtain the six-cycle $\left\{\vectr{0}{1},\vectr{1}{-1},\vectr{-1}{-1},\vectr{0}{-1},\vectr{-1}{1},\vectr{1}{1}\right\}$. This cycle corresponds to a different periodic word:   $$\underline{\vectr{3}0\vectr0{-2}\vectr30\vectr{-3}0\vectr0{2}\vectr{-3}{0}}.$$
This means that the 6-cycle is traveled by two diferent paths according to the matrix $A$ or $A^T$. 
\end{example}

\subsection{Tiling and spectra in some examples}
\par One of the uses of encoding is applications to tiling questions in $\br^d$. The simplest tiles $X$ in $\br^d$ are measurable subsets which make a tiling of $\br^d$ by translations with vectors from some lattice, say $\Gamma$ (i.e, a rank d subgroup). Since we work in the measurable category we allow different translates in the tiling $X + \gamma$, for $\gamma\in\Gamma$ to overlap on sets of measure zero.

One might think that when $X$ is given, then the presence of a suitable lattice $\Gamma$ making $X$ into a translation tile for $\br^d$ could be decided by visual inspection, at least in the case of $d = 2$.  After all, when a pair $(A,\mathcal D)$ is given, then there are fast Mathematica programs which produce excellent plots of the attractor sets $X = X(A^T,\mathcal D)$, black on white; see for example \cite{BrJo99}. But except for isolated cases, it turns out that when the planar sets $X$ are represented in black on white, then there will typically be many white spots, or gaps, disconnecting $X$ in complicated ways. If some lattice $\Gamma$ will make $X$ into a translation tile, then the white areas must be filled in by black under translations $X + \gamma$, $\gamma \in \Gamma$. An inspection of \cite{BrJo99} reveals that this is not easy to discern by visual inspection. Hence, instead we resort below  to spectral theoretic tools for locating the lattices which do the job. 

\par A more complicated form of tilings (still with a single base tile) refer to the case when the set $\Gamma$ of translation vectors is a set which is not a lattice, e.g., translation sets of quasiperiodic tilings. But for our present considerations lattice tilings will suffice. 
\par The sets X which will interest us are the attractors $X = X(A^T,\mathcal D)$ from affine IFSs as described in Corollary \ref{corsum}. It is known that every such $X$ is compact with non-empty interior, and so in particular it has positive d-dimensional Lebesgue measure. 
\par Hence it is of interest to ask for a spectral analysis of the Hilbert space $L^2(X)$, referring to $d$-dimensional Lebesgue measure. In fact, using Pontryagin duality for abelian groups, one can check that $X$ tiles $\br^d$ with a lattice if and only if the dual lattice makes an orthogonal basis of complex exponentials in $L^2(X)$. The result is often refered to as Fuglede's theorem. For background, see the references \cite{Fug74} and \cite{Rud62}.

\par To understand the correspondence between choice of translation lattice on the one hand and spectrum on the other we need:

\begin{definition}
Let $X\subset\br^d$ be measurable with $0<\mu(X)<\infty$ where $\mu$ denotes the $d$-dimensional Lebesgue measure. For $\xi\in\br^d$ set 
$e_{\xi}(x)=e^{2\pi i\xi\cdot x}$, where $\xi\cdot x:=<\xi\mbox{, }x>=\xi_1x_1+\cdots \xi_dx_d$ and $x=(x_1,\cdots,x_d)\in\br^d$. If $\Gamma\subset\br^d$ is a discrete subgroup (in this case a rank $d$ lattice) set 
$$E_X(\Gamma):=\{ e_{\xi}|_X\mbox{ : }\xi\in\Gamma\}$$ where $|_X$ denotes restriction to the set $X$. \\
If $\Gamma$ is a lattice we set 
$$\Gamma^{\circ}:=\{\lambda\in\br^d\mbox{ }|\mbox{ }\lambda\cdot\xi\in\bz\mbox{ for all }\xi\in\Gamma\}$$
called the dual lattice. \\
If $\Lambda\subset\br^d$ is a discrete subset we say that it is a {\it spectrum} for $X$ or that the pair $(X, \Lambda)$ is a {\it spectral pair}
 iff $E_X(\Lambda)$ is an orthogonal basis in the Hilbert space $L^2(X)=L^2(X,\mu)$.
\end{definition}

\begin{lemma}\label{lemfuglede} (Fuglede \cite{Fug74}) Let $0<\mu(X)<\infty$ and $\Lambda$ a rank-$d$ lattice. The following conditions are equivalent:\\
i) $E_X(\Lambda)$ is an orthonormal basis in $L^2(X)$;\\
ii) $X$ tiles $\br^d$ by the dual lattice $\Lambda^{\circ}$.
\end{lemma}

\begin{remark}
We can draw the following stronger conclusion: When $X$ is given, there are no other tiling lattices for $X$ than those which arise as in (ii) by spectral duality.
       Proof. Every lattice $\Lambda$ satisfies ${\Lambda^\circ} ^{\circ} = \Lambda$, i.e., the double dual yields back the initial lattice. To see this, use the following general observations which also serve to make explicit the standard lattice operations which we will be using in the proof. 
\end{remark}
Referring to the IFS of Definition \ref{ifs} we note the following formula for the computation of the $L^2(X(A^T,\mathcal D))$-inner products. Set $X=X(A^T,\mathcal D)$ and 
$\hat \chi_X(\xi):=\int_Xe_{\xi}(x)dx$. Then 
$$\hat \chi_X(\xi)=\prod_{n=1}^{\infty}m_{\mathcal D}(A^{-n}\xi)$$
where $$m_{\mathcal D}(\xi):=\frac{1}{|\det A|}\sum_{d\in\mathcal D}e_d(\xi).$$ Recall that $|\det A|=$number of elements in $\mathcal D$. 
\par Some remarks about lattices in $\br^d$ are in order: every lattice is by definition a rank-$d$ subgroup of $\br^d$ and it can be shown that it has the form $\Gamma=M\bz^d$, where $M$ is an invertible $d\times d$ matrix and where points in $\bz^d$ are represented by column vectors. We will write $\Gamma_M$ to emphasize the matrix $M$ that completely determines the lattice. The next lemma is elementary:

\begin{lemma}
{\rm (i)} $\Gamma_M\subset\bz^d$ if and only if $M\in\mathcal M_d(\bz)$.

{\rm (ii)} If $\Gamma=\Gamma_M$ then $\Gamma^{\circ}=\Gamma_{(M^T)^{-1}}$. In other words if $\Gamma$ is given by $M$ then its dual is 
given by $(M^T)^{-1}$.

{\rm (iii)} $\Gamma^{\circ\circ}=\Gamma$
\end{lemma}

\par We will use names from \cite{BrJo99} for the fractals $X = X(A^T,\mathcal D)$ in $\br^2$. 
These names refer both to their geometric appearance as planar sets $X$, as
well as to a counting of $\bz^2$-cycles, i.e., the number of points in 
$(-X)\cap\bz^2$. See \cite{BrJo99} (end of subsection 9.3) for details. For example, Cloud-Nine has
three one-cycles and one six-cycle in $\bz^2$.

What follows is a family of examples in 2D. In each case, we are
asking the following questions: How much flexibility is there in
selecting digits when the base for the 2D number system is fixed? In
our case, we are using the positional radix representation for
vectors, and thus the base for our number system is a chosen matrix $A$.
For several of the examples below, we fix a particular $A$, and then we
vary our choices of ``digit'' sets $\mathcal D$ in $\bz^2$. The points in $\mathcal D$ will serve
as ``digits'' in a positional representation.

   We are motivated by Knuth's algorithmic approach mentioned in the
Introduction: What are the ``integers'' and what are the ``fractions'' in
a number system specified by a particular pair $(A,\mathcal D)$? What is the
encoding, and what is the decoding?
   When the matrix $A$ is fixed, how do changes in $\mathcal D$ reflect themselves
in the answer to the questions?

   The examples below are sketched with Mathematica programming in
\cite{BrJo99}, and the names we use for the fractals $X$ are consistent with
\cite{BrJo99}, i.e., the the respective names of the sets $X$, Cloud-Nine etc.
The examples when $A$ is the same but $\mathcal D$ changes are referred to by the
name Cloud, followed by a number. The number indicates the cardinality
of $(-X)\cap \bz^2$.

    However the questions addressed here are different from those of [BJ99].

    It is of interest to understand how much flexibility there is in
selecting digits when the base for the number system is fixed. In our
case, the base for our vector number system is the matrix $A$, and so we
vary the choices for the companion set $\mathcal D$. But when $A$ is given, the
choice of $\mathcal D$ is always restricted by demanding a bijection $\mathcal D \leftrightarrow
\bz^2/A^T\bz^2$.

\begin{itemize}
    \item  Cloud Three. $A=\left(\begin{array}{cc}1&-2\\2&1\end{array}\right)$, $\mathcal D=\left\{\vectr00,\vectr0{\pm 1},\vectr{0}{\pm2}\right\}$. Here the lattice $2\bz\times\bz$ 
makes $X$ tile $\br^2$. Cloud Three only has one-cycles on $\bz^2$, i.e., $(-X)\cap\bz^2= \mathcal C_1$. Moreover $X$ is not a Haar wavelet. It has measure = 2.

    \item  Cloud Five. $A=\left(\begin{array}{cc}1&-2\\2&1\end{array}\right)$, $\mathcal D=\left\{\vectr00,\vectr{\pm3}0,\vectr{\pm1}0\right\}$. Cloud Five is a lattice tile with
lattice $\bz\times 2\bz$. So
$X$ is not a Haar wavelet. It has measure = 2.

     \item Cloud Nine. $A=\left(\begin{array}{cc}1&-2\\2&1\end{array}\right)$, $\mathcal D=\left\{\vectr00,\vectr{\pm3}{0},\vectr{0}{\pm2}\right\}$. Cloud Nine is a lattice tile with the lattice $\bz\times 2\bz$.  Cloud Nine $X$ has three one-cycles and one six-cycle. So Cloud Nine is not a Haar wavelet. It has measure = 2.

      \item Twin Dragon. $A=\left(\begin{array}{cc}1&1\\-1&1\end{array}\right)$, $\mathcal D=\left\{\vectr00,\vectr10\right\}$. The Twin Dragon is a lattice tile with lattice $\bz^2 = \Gamma$. So the
Twin Dragon is a Haar wavelet. It has measure = 1.
\end{itemize}

 We are using Lemma \ref{lemfuglede} in identifying lattices which make the
various Cloud examples $X$ tile $\br^2$. For this purpose we must identify
our cycles relative to so called Hadamard systems as defined in
\cite{DuJo07c}. A Hadamard system consists of a matrix $A$ and two sets $\mathcal D$ and $\mathcal L$
as dual digits, $\#\mathcal D = \#\mathcal L = |\det A|$. By ``dual'' we mean that the matrix
formed from the exponentials as
\begin{equation}\label{eqhada}
\frac1{\sqrt{|\det A|}} (e^{ 2 \pi i \at^{-1}d\cdot l})_{ d \in\mathcal D, l \in\mathcal L}.
\end{equation}
is a unitary $|\det A|\times |\det A|$ matrix.

Let $N$ be the absolute value of the determinant, and let $\bz_N$ be the cyclic group of order $N$. Then the matrix $U_N$ for the Fourier transfrom on $\bz_N$ is an example of a Hadamard matrix as in \eqref{eqhada}; specifically the $j,k$ entry in $U_N$ is  $\frac{1}{\sqrt{N}} \zeta^{jk}$, $j,k \in \bz_N$, where $\zeta= \zeta_N$  is a fixed principal $N$'th root of $1$.
\begin{proof}
To find the lattices that give tiles for these examples, we use Lemma \ref{lemfuglede}, and find the dual lattices that give orthogonal bases of exponentials. For this we use the techniques introduced in \cite{DuJo06,DuJo07c}. 

First let us look at the matrix $A=\left(\begin{array}{cc}1&-2\\2&1\end{array}\right)$ for the cloud examples. We want to find what is the lattice that makes $X(A^T,\mathcal D)=\{\sum_{j=1}^\infty\at^{-j}d_j\,|\,d_j\in\mathcal D\}$ tile $\br^2$. For this we need a set $\mathcal L$ such that $\frac{1}{\sqrt{|\det A|}}(e^{2\pi i \at^{-1}d\cdot l})_{d\in\mathcal D,l\in\mathcal L}$ is a unitary matrix, i.e., $(A^T, \mathcal D,\mathcal L)$ is a Hadamard triple. It is enough to take $\mathcal L$ a complete set of representatives for $\bz^2/A\bz^2$. We will take $\mathcal L$ to be 
$$\mathcal L:=\left\{\vectr00,\vectr{\pm3}{0},\vectr{0}{\pm2}\right\}.$$

With this choice of $\mathcal L$ the reader may check that for each of the Cloud-examples listed above, the corresponding Hadamard matrix from \eqref{eqhada} turns out, up to permutation, to simply agree with the matrix $U_5$ for the Fourier transform on $\bz_5$.

According to \cite{DuJo07c} we have to see if there are any proper invariant subspaces for $A$. But those would give rise to real eigenvalues of $A$, and this is not the case. Thus, by \cite{DuJo07c}, the spectrum of $X(A^T,\mathcal D)$, is determined only by the ``$m_{\mathcal D}$-cycles''. These are the cycles $C=\{x_0,\dots,x_{p-1}\}$ for the ``dual'' IFS $\sigma_l(x)=A^{-1}(x+l)$, $l\in \mathcal L$ with the property that $|m_{\mathcal D}(x_i)|=1$ for all $i\in\{0,\dots,p-1\}$.

Then, by \cite{DuJo06,DuJo07c} the spectrum of $X(A^T,\mathcal D)$ is the smallest set $\Lambda$ that contains $-C$ for all the $m_{\mathcal D}$-cycles, and such that $A\Lambda+\mathcal L\subset \Lambda$.

{\bf Cloud Three}: We have 
$$m_{\mathcal D}(x,y)=\frac15\left(1+e^{2\pi iy}+e^{-2\pi iy}+e^{2\pi i2y}+e^{-2\pi i2y}\right).$$
If we want $|m_{\mathcal D}(x,y)|=1$ then we must have that all the terms in the sum are $1$ so $y\in\bz$, and $x$ is arbitrary. 

We are looking for $m_{\mathcal D}$-cycles, so $|m_{\mathcal D}(\sigma_l(x,y))|=1$ for some $l\in\mathcal L$ so $\sigma_l(x,y)$ must have the second component in $\bz$. The inverse of $A$ is 
$A^{-1}=\frac{1}{5}\left(\begin{array}{cc}1&-2\\2&1\end{array}\right)$. Thus $\frac15(2(x+l_x)+(y+l_y))\in\bz$. This implies that $x\in\frac12\bz$. 

We claim that $\Lambda=\frac12\bz\times\bz$. For this, note first that $A(\frac12\bz\times\bz)+\mathcal L\subset \frac12\bz\times\bz$. By the previous computation, $\frac12\bz\times\bz$ contains the negative of all the $m_{\mathcal D}$-cycles. Then, take $w_0:=(\frac k2,k')\in\frac12\bz\times\bz$. Then for $l\in\mathcal L$,
$$\sigma_l(-\frac k2,-k')=(\frac1{10}(-k+2l_x+4k'-4l_y),\frac15(-k+2l_x-k'+l_y)).$$
Note that there is a unique $l_0\in\mathcal L$ such that $w_1:=-\sigma_{l_0}(-\frac k2,-k')\in\frac12\bz\times\bz$. As in the proof of Proposition \ref{proprc}, there is a sequence $l_0,l_1\dots\in\mathcal L$ such that if $w_{n+1}=-\tau_{l_n}(-w_n)$, $n\in\bn$, then $w_n\in\frac12\bz\times\bz$ and, for some $m$, $-w_m$ is a cycle point for $(\sigma_l)_{l\in\mathcal L}$. Note that since $w_n\in\frac12\bz\times\bz$, $-w_m$ is a point in a $m_{\mathcal D}$-cycle. 

Since $w_{m}=-\tau_{l_{m-1}}(w_{m-1})$, we have that $w_{m-1}=Aw_{m}+l_{m-1}\in A(-C)+\mathcal L$, where $C$ is the $m_{\mathcal D}$-cycle of $-w_m$. By induction we obtain that $w_{0}$ must be in $\Lambda$. Thus $\Lambda=\frac12\bz\times\bz$ is the spectrum. Taking the dual we obtain that $X(A^T,\mathcal D)$ tiles $\br^2$ by $2\bz\times\bz$.

{\bf Cloud Five}: We have

$$m_{\mathcal D}(x,y)=\frac15\left(1+e^{2\pi i3x}+e^{-2\pi i3x}+e^{2\pi ix}+e^{-2\pi ix}\right).$$
Therefore $|m_{\mathcal D}(x,y)|=1$ iff $x\in\bz$. For $m_{\mathcal D}$-cycles we must have that the first component $\sigma_l(x,y)$ must be in $\bz$ for some $l\in\mathcal L$. This implies that $y\in\frac12\bz$.

We claim that $\Lambda=\frac12\bz\times\bz$. The proof works just as for the Cloud Three example so we will leave it to the reader. Thus the dual lattice is $2\bz\times\bz$ and $X(A^T,\mathcal D)$.

{\bf Cloud Nine}: We have
$$m_{\mathcal D}(x,y)=\frac15\left(1+e^{2\pi i3x}+e^{-2\pi i3x}+e^{2\pi i2y}+e^{-2\pi i2y}\right).$$
So $|m_{\mathcal D}(x,y)|=1$ iff $x\in\frac13\bz$ and $y\in\frac12\bz$. For $m_{\mathcal D}$-cycles we must have $(x,y)=(\frac k3,\frac {k'}2)$ with $k,k'\in\bz$ and also $|m_{\mathcal D}(\sigma_l(\frac k3,\frac{k'}2))|=1$ for some $l\in\mathcal L$. This implies that $\frac15(2(\frac k3+d_x)-(\frac{k'}2+d_y))\in\frac12\bz$ so $k$ must be divisible by $3$. Thus the $m_{\mathcal D}$-cycles are contained in $\bz\times\frac12\bz$. 

Just as in the previous examples we get that $\Lambda=\bz\times\frac12\bz$ so the tiling lattice is $\bz\times 2\bz$.

{\bf Twin Dragon}: For $A=\left(\begin{array}{cc}1&1\\-1&1\end{array}\right)$, there are no proper invariant subspaces so the analysis of the $m_{\mathcal D}$-cycles will suffice. We can take $\mathcal L:=\left\{\vectr00,\vectr10\right\}$.

$$m_{\mathcal D}(x,y)=\frac12(1+e^{2\pi ix}).$$
Therefore $|m_{\mathcal D}(x,y)|=1$ iff $x\in\bz$. $A^{-1}=\frac12\left(\begin{array}{cc}1&-1\\1&1\end{array}\right)$. We want 
the first component of $\sigma_l(x,y)$ to be in $\bz$ so $\frac12(x-y)\in\bz$, therefore $y\in\bz$. 

As in the previous examples we can check that $\Lambda=\bz^2$ so the dual lattice is $\bz^2$.
\end{proof}

The method used in our analysis of the examples may be formalized as follows. Stated in general terms it applies to a large class of IFSs which carry a fairly minimal amount of intrinsic duality. For the convenience of the reader, we begin with two definitions. This is of interest as there are few results in the literature which produce formulas for lattices which turn particular attractors $X$ into tiles under the corresponding translations. For details we refer to \cite{DuJo07c,DuJo06}.

\begin{definition} (i) A Hadamard triple in $\br^d$ is a system $(A,\mathcal D, \mathcal L)$ where $A$ is expansive in $\mathcal M_d(\bz)$, $\mathcal D,\mathcal L$ are in $\bz^d$, and $\mathcal L$ is such that the matrix \eqref{eqhada} is unitary.

(ii) For a Hadamard triple $(A,\mathcal D,\mathcal L)$ the cycles $C$, that correspond to the IFS $(\sigma_l)_{l\in\mathcal L}$, for which the absolute value of $m_{\mathcal D}$ is 1 are called extreme relative to $m_{\mathcal D}$, or $m_{\mathcal D}$-cycles.
\end{definition}

We are now ready to state our general tiling result.

\begin{corollary}
Let $(A,\mathcal D,\mathcal L)$ be a Hadamard triple in $\br^d$, and such that $\mathcal D$ is a complete set of representatives for $\bz^d/A^T\bz^d$, and let $X = X(A^T,\mathcal D)$ be the corresponding Hutchinson attractor. Suppose $A$ has no proper invariant subspaces.
Let $\Lambda$ be the smallest lattice in $\br^d$ containing all the sets $-C$ where $C$ runs over the $m_{\mathcal D}$-extreme cycles, and which is invariant under the affine mappings $x\mapsto Ax + l$, for $l \in \mathcal L$.
Then the dual lattice $\Gamma = \Lambda^\circ$ makes $X$ tile $\br^d$ with $\Gamma$ translations.
\end{corollary}
\begin{acknowledgements}
The second named author had helpful discussions with Prof Sergei Silvestrov, University of Lund, Sweden. Useful suggestions from a referee led to
improvements in the presentation.
\end{acknowledgements}

\bibliographystyle{alpha}
\bibliography{encod}
\end{document}